\documentclass[11pt]{article}
\usepackage{amsmath,amssymb,amsthm,mathrsfs}
\usepackage{enumitem}
\usepackage{geometry}
\usepackage{xcolor}
\usepackage{url}
\numberwithin{equation}{section}

\newtheorem{theorem}{Theorem}[section]
\newtheorem{lemma}[theorem]{Lemma}
\newtheorem{proposition}[theorem]{Proposition}
\newtheorem{conjecture}[theorem]{Conjecture}
\newtheorem{corollary}[theorem]{Corollary}
\newtheorem{remark}[theorem]{Remark}

\newtheorem{problem}[theorem]{Problem}

\newcommand{\E}{\mathrm{e}}
\newcommand{\D}{\mathrm{d}}

\newcommand{\vertk}{\stackrel{{\cal D}}{\longrightarrow}}
\newcommand{\stk}{\stackrel{{\mathbb{P}}}{\longrightarrow}}
\newcommand{\edist}{\stackrel{{\cal D}}{=}}
\title{Extreme-value asymptotics for smallest enclosing balls of spherically symmetric samples}
\author{Norbert Henze\footnote{Institute of Stochastics, Karlsruhe Institute of Technology (KIT), 76131 Karlsruhe, Germany. e-mail:
henze@kit.edu} \ and Michael Klatt\footnote{German Aerospace Center (DLR), Institute for AI Safety and Security,
Wilhelm-Runge-Str. 10, 89081 Ulm, Germany;
German Aerospace Center (DLR), Institute of Frontier Materials on Earth and in Space,
Functional, Granular, and Composite Materials, 51170 Cologne, Germany;
Department of Physics, LMU München, Schellingstr. 4, 80799 Munich, Germany. e-mail: Michael.klatt@dlr.de}}

\date{}

\begin{document}
\maketitle

\begin{abstract}
Let $R_n$ denote the radius of the smallest ball containing an i.i.d.\ sample of spherically symmetric random vectors in $\mathbb R^d$, and let $D_n$ be the sample diameter. While the asymptotic behaviour of $D_n$ is well understood, much less is known about $R_n$. Starting from a classical result of Daniels for planar Gaussian samples, we develop a point-process approach to the asymptotic behaviour of $R_n$ in the three classical extreme-value regimes. We obtain limit theorems in the Gumbel, Fr'echet and Weibull cases, with the limiting laws described by geometric functionals of Poisson point processes. In dimension two, the Gumbel limit can be identified explicitly and is universal within the corresponding max-domain of attraction. We also study the asymptotic relationship between $R_n$ and $D_n$ and the supporting configurations of the smallest enclosing ball.
\end{abstract}




\section{Introduction}

Let $X, X_1, X_2, \ldots$ be a sequence of independent and identically distributed (i.i.d.) $d$-dimensional
random vectors (“points”) with an absolutely continuous distribution with respect to Lebesgue measure.
Let $\|\cdot\|$ denote the Euclidean norm on $\mathbb{R}^d$, and let
$
B(x,r) = \{y \in \mathbb{R}^d : \|x-y\| \le r\}
$
denote the closed ball centred at $x \in \mathbb{R}^d$ with radius $r>0$. The enclosing ball of $X_1,\ldots,X_n$
is the smallest ball $B(x,r)$ that contains $\{X_1,\ldots,X_n\}$. Writing
\[
R_n := \inf_{c\in\mathbb{R}^d} \max_{1\le j\le n} \|X_j-c\|
\]
for the radius of this ball, a natural question is whether $R_n$, suitably normalized,
has a nondegenerate limit distribution as $n\to\infty$, and, if so, how this limit law is related to the geometric structure of the sample.
The statistic $R_n$  arises naturally in computational geometry and stochastic geometry,
where the smallest enclosing ball provides a basic descriptor of the spatial extent of a point cloud.

If $d=1$, then $R_n$ equals half the sample range
$
\max_{1\le j\le n} X_j - \min_{1\le j\le n} X_j$.
Since the range has been studied extensively (see, e.g., \cite{dav81}, Section 9.4), we restrict attention to the case $d\ge2$ throughout.

A natural quantity closely related to $R_n$ is the diameter (maximum interpoint distance)
\[
D_n = \max_{1\le i<j\le n} \|X_i - X_j\|.
\]
The relation between these two random variables is governed by the classical inequality
\begin{equation}\label{jung}
R_n \le D_n \sqrt{\frac{d}{2(d+1)}}
\end{equation}
due to Jung \cite{jung01}.
Equality in \eqref{jung} holds if and only if the enclosing ball is determined by $d+1$ boundary points
forming the vertices of a regular simplex with edge length $D_n$. In particular, for $d=2$ and $d=3$, equality occurs precisely
when the defining points form an equilateral triangle or a regular tetrahedron with edge length $D_n$.
Under absolute continuity, such configurations have probability zero.

On the other hand, we trivially have
\begin{equation}\label{dn-less-rn}
\frac{D_n}{2} \le R_n.
\end{equation}
Together, \eqref{jung} and \eqref{dn-less-rn} show that $R_n$ and $D_n$ are of the same order as $n\to\infty$.
Moreover, the points lying on the boundary of the enclosing ball necessarily belong to the set of vertices of
the convex hull $\mathcal{C}_n$ of ${X_1,\ldots,X_n}$. Hence $R_n$ is a functional of $\mathcal{C}_n$. For an overview of functionals
of random convex hulls we refer to Reitzner \cite{rei10}.

While the asymptotic behaviour of the diameter $D_n$ is well understood
(see, e.g., \cite{anr}, \cite{dfs15}, \cite{he26}, \cite{heja}, \cite{henzelao}, \cite{hekl96}, \cite{jamsva15}, \cite{lao10}, \cite{mr93},
\cite{maymol07}, \cite{sc16} and \cite{sc19}),  corresponding results for the enclosing radius $R_n$ appear to be largely unavailable.  First-order asymptotics can, however, be deduced in some cases from general results on random polytopes. For instance,
Theorem~2.3 of Fresen and Vitale \cite{frvi} implies for Gaussian samples that
$R_n/\sqrt{2\log n}\rightarrow 1$ almost surely. 
The only explicit distributional result appears to be a classical
formula due to Daniels \cite{dan52} for the planar Gaussian case, see Section~\ref{secdaniels}.
From this formula one obtains a non-Gumbel limit law for suitably normalized
$R_n$. By contrast, the diameter
$D_n$  has a Gumbel limit distribution under the same assumptions
(Jammalamadaka and Janson \cite{jamsva15}). This discrepancy already indicates that the
statistics $2R_n$ and $D_n$ cannot in general share the same asymptotic distribution.
Given the apparent difficulty of deriving explicit limit laws for
$R_n$, the present work develops a general asymptotic framework and
also formulates conjectures and open problems that may guide future
research. 

The present paper investigates the asymptotic behaviour of the enclosing radius $R_n$ for samples
from spherically symmetric distributions.
Recall that the distribution of $X$ is called \emph{spherically symmetric} if it admits the  representation $X \edist \rho U$,
where $\edist$ denotes equality in distribution, $\rho =\|X\|$, $U$ is uniformly distributed
on the unit sphere $\mathbb{S}^{d-1}= \{x \in \mathbb{R}^d: \|x\|=1\}$, and $\rho$ and $U$ are independent.

The remainder of the paper is organized as follows. Section~\ref{secdaniels} revisits Daniels' classical result for the radius of the enclosing circle in
the planar Gaussian case and derives its asymptotic consequences.
Section~\ref{secgumbelkotz} develops a general asymptotic framework for the enclosing radius $R_n$ in the Gumbel regime. When  the radial
distribution belongs to the Gumbel max-domain of attraction, we derive
a Poisson functional limit theorem for $R_n$ based on a geometric
linearization and a variational representation involving a limiting
Poisson point process. In dimension two, the limiting distribution is
shown to be universal throughout the Gumbel max-domain of attraction
and coincides with the limit law obtained from Daniels' classical
Gaussian formula. We also show explicitly that the Kotz family,
including the Gaussian distribution, satisfies the assumptions of the
general theorem.

Section~\ref{secfrechet} studies the case of heavy-tailed radial distributions and derives an asymptotic
result for $R_n$ in the Fréchet regime, complementing the corresponding limit theory for $D_n$,
which was obtained independently by Jammalamadaka and Janson \cite{jamsva15} and by Henze and Lao \cite{henzelao}.

In Section~\ref{secweib}, we study the enclosing radius in the bounded-support (Weibull) regime. Assuming polynomial decay of the radial distribution near the boundary of the  unit ball, we establish convergence of a boundary point process and combine a geometric
linearization with a variational representation to derive a Poisson functional limit theorem for suitably normalized versions of $1-R_n$. As a special case, we obtain the limit distribution of $n(1-R_n)$ for samples uniformly distributed on the unit ball. We
also discuss scaling properties of the limiting functional, compare the normalization with that of the diameter, and provide a numerical illustration of the limit theorem.

Section~\ref{secratio} investigates the relationship between the enclosing radius $R_n$ and the diameter $D_n$,
establishing ratio theorems.

Section~\ref{secsuffcond}  studies the asymptotic behaviour of the number of supporting points of the
smallest enclosing ball. In the planar case we show that the three-point configuration is
asymptotically dominant under suitable conditions on the radial distribution, and we establish the same
behaviour for the Gaussian model. In higher dimensions, we provide simulation evidence
 suggesting dominance of the $(d+1)$-point configuration in the bounded-support setting. We conclude with
a conjecture and an open problem concerning the Gumbel case.

Section~\ref{secproblems} collects the main open problems that remain
after the preceding results. These concern, among other questions,
more explicit descriptions of the limiting Poisson functionals, the
asymptotic behaviour of the centres of the smallest enclosing balls,
supporting configurations in higher dimensions, finer comparisons
between the enclosing radius and the diameter, and extensions beyond
spherical symmetry.

Throughout the paper, the symbols $\vertk$ and $\stk$
denote convergence in distribution and convergence in probability, respectively.


\section{Daniels' result and its asymptotic implications}\label{secdaniels}


Daniels \cite{dan52} considered i.i.d.\ random points from a bivariate spherically symmetric normal distribution.
He allowed for a general variance of the components, but we assume unit variance in what follows, as the general case reduces to this by scaling.
By distinguishing the two main categories that (i) the enclosing circle
$B_n$, say, passes through three of the points forming an acute-angled triangle,
and the remainder of the points lie inside the circle, and (ii) the circle passes through two
diametrically opposite points, the remainder lying inside the circle, and by studying the infinitesimal probability that the
radius of $B_n$ lies between $r$ and $r + \textrm{d} r$ and that the centre of $B_n$ is at distance
$c$ to $c + \textrm{d} c$ from the origin, \cite{dan52} obtained the joint density of $R_n$
and the distance of the centre of $B_n$ from the origin. From this result, the distribution function
$F_n(r) = \mathbb{P}(R_n \le r)$ was obtained, and is given by
\begin{equation}\label{distrn}
F_n(r) \ = \ n \left(1- \exp\left(- \frac{r^2}{2} \right)\right)^{n-1} - (n-1)\left(1- \exp\left(- \frac{r^2}{2} \right) \right)^n, \quad 0 \le r < \infty.
\end{equation}
This explicit formula appears to be the only available closed-form distributional result
for the enclosing radius $R_n$, and thus provides a natural benchmark for the asymptotic
analysis developed in the present paper. It highlights the exceptional role of Daniels'
result in a setting where explicit distributional information is, to the best of our
knowledge, largely unavailable.

Interestingly, as also noticed in \cite{dan52}, $F_n$ is the distribution function of the radius of the smallest circle,
{\em centred at the origin}, enclosing {\em all but one} of $X_1,\ldots,X_n$.

Indeed, denoting by $M_{n,2}$ the {\em second largest} value among $\|X_1\|, \ldots,\|X_n\|$, and writing
$S_n = \sum_{j=1}^n {\bf 1}\{\|X_j\| >r \}$ 
for the number of points that lie outside a circle centred at the origin with radius $r$, we have $S_n \le 1$ if, and only if, $M_{n,2} \le r$.
Since the distribution of $S_n$ is binomial with parameters $n$ and $p= \mathbb{P}(\|X\| >r)$, and since the distribution of $\|X\|^2$ is exponential
with expectation $2$, this result follows readily.

The explicit formula \eqref{distrn} provides a convenient starting point for studying the asymptotic behaviour
of $R_n$. As noted above, $R_n$ coincides in distribution with the second largest order statistic of
$\|X_1\|,\ldots,\|X_n\|$. This representation naturally connects $R_n$ to extreme-value theory. In particular,
one expects that the limit behaviour of $R_n$ can be derived from the well-known asymptotics of the maximum norm.
The following result makes this precise.

\medskip

\begin{theorem}[Limit law of $R_n$ in the planar standard normal case]\label{haupt}
Let $X_1,X_2, \ldots $ be i.i.d. points in $\mathbb{R}^2$ having a centred normal distribution with unit covariance matrix.
We then have
\[
\lim_{n \to \infty} \mathbb{P} \left(\sqrt{2 \log n} \left(R_n - \sqrt{2 \log n} \right) \le t \right) =  G(t)(1+ \exp(-t)), \quad t \in \mathbb{R},
\]
where $G(t) = \exp(-\exp(-t))$, $t \in \mathbb{R}$, is the distribution function of the Gumbel distribution.
\end{theorem}    

\medskip
\noindent The proof of Theorem~\ref{haupt} uses the following result on the maximum norm $M_n \ = \ \max_{j=1,\ldots,n} \|X_j \|$ of $X_1,\ldots, X_n$, which holds in any dimension $d \ge 2$.

\begin{proposition}[Limit law of $M_n$ in the spherically symmetric normal case]\label{propo1}
Let $X, X_1,X_2, \ldots $ be i.i.d. random points in $\mathbb{R}^d$ having a centred normal distribution with unit covariance matrix.
Then
\[
\lim_{n \to \infty} \mathbb{P} \left(\sqrt{2 \log n}\left(M_n - a_n \right) \le t \right) \ = \ G(t), \qquad t \in \mathbb{R},
\]
where
\[
a_n = \sqrt{2 \log n} + \frac{(d-2)\log \log n}{2 \sqrt{2 \log n}} - \frac{\log \Gamma(d/2)}{\sqrt{2 \log n} }.
\]
\end{proposition}

\begin{proof}
Since $\|X\|^2$ has a $\chi^2_d$-distribution and thus a Gamma distribution, the distribution of $M_n^2$
is in the maximum domain of attraction of the Gumbel law. From p. 156 of~\cite{ekm97} we therefore obtain
$
(M_n^2 - b_n)/2 \ \vertk \ Z,
$
where
\[
b_n \ = \ 2 \left( \log n + \left( \frac{d}{2} - 1 \right) \log \log n - \log \Gamma(d/2) \right),
\]
and $Z$ has a Gumbel distribution.  Using $M_n^2- b_n = (M_n - \sqrt{b_n})(M_n + \sqrt{b_n})$ together with
$(M_n + \sqrt{b_n})/\sqrt{b_n} \stk 2$ and Slutsky's lemma, we have $\sqrt{b_n}(M_n - \sqrt{b_n}) \vertk Z$.
Upon noting that $\sqrt{b_n} = a_n + o(1/\sqrt{2 \log n})$ as $n \to \infty$, the assertion follows.
\end{proof}

\medskip

\noindent {\em Proof} of Theorem~\ref{haupt}. Notice that for the special case $d=2$ Proposition  \ref{propo1} reads
\[
\sqrt{2 \log n} \left(M_n - \sqrt{2 \log n}\right) \ \vertk \ Z,
\]
where $Z$ has a Gumbel distribution. Since, as noted above, $R_n$ has the same distribution as the second largest
of the order statistics of $\|X_1\|, \ldots, \|X_n\|$, the result follows  from a standard argument
of extreme value theory (see, e.g. \cite{ekm97}, p. 198).
\qed

\begin{remark}\label{remdiamg} For the diameter $D_n$ of a Gaussian sample, the situation is different. From  Corollary 1 of \cite{mr93} we have
\[
\lim_{n \to \infty} \mathbb{P}\left(2\sqrt{2 \log n} \left( \frac{D_n}{2} - \sqrt{2 \log n}
+ \frac{\ell_2 n - 2 \ell_3n + 5 \log 2 + \log \pi}{8 \sqrt{2 \log n}} \right) \le t \right) = G(t),
\]
where $\ell_2 n = \log \log n $ and $\ell_3  n = \log \log \log n$.

More generally, the Gumbel limit behaviour of the maximum interpoint distance for samples from
a spherically symmetric distribution was established by Jammalamadaka and Janson \cite{jamsva15}. In contrast to this
classical extreme-value behaviour of the diameter, Theorem~\ref{haupt} shows that the enclosing radius
$R_n$ obeys a different limit law even in the Gaussian case. In particular, the lower bound $D_n/2$ is asymptotically
slightly smaller than $R_n$.
\end{remark}

\begin{remark}\label{remcombine}
Combining Theorem~\ref{haupt} with Remark~\ref{remdiamg}, we obtain
\[
R_n-\frac{D_n}{2}\stk 0 .
\]
Thus, although $R_n$ and $D_n$ have different asymptotic limit laws, their difference is asymptotically negligible.
In particular, this implies that
\[
\frac{R_n}{D_n/2} \stk 1.
\]
We return to this ratio in Section~\ref{secratio}, where its asymptotic behaviour
is studied under more general assumptions on the underlying distribution.

It is an open problem whether there exist sequences $a_n \to \infty$ and $b_n \ge 0$ such that
\[
a_n\left(R_n - \frac{D_n}{2} - b_n\right)
\]
has a nondegenerate limit distribution as $n \to \infty$.
In particular, it is unclear whether a canonical normalization exists.
\end{remark}

%

\section{The Gumbel case}\label{secgumbelkotz}
%
%

In this section we study the radius $R_n$ of the smallest enclosing ball
when the radial distribution belongs to the Gumbel max-domain of
attraction. The arguments are based on point-process convergence and a
geometric linearization of the enclosing-ball problem.

Let $X,X_1,X_2,\ldots$ be i.i.d.\ spherically symmetric random vectors
in $\mathbb{R}^d$, $d\geq 2$, and write
$X \edist \rho U$,
where $\rho=\|X\|$, $U$ is uniformly distributed on
$\mathbb{S}^{d-1}$, and $\rho$ and $U$ are independent. Accordingly,
we write $X_i=\rho_i U_i$, $i\geq 1$.

We assume that there exist sequences $a_n>0$ and $b_n>0$ such that
\begin{equation}\label{eq:gumbel_tail}
 n\,\mathbb{P}\bigl(\rho>a_n+b_n t\bigr)
 \rightarrow \E^{-t},
 \qquad t\in\mathbb{R}.
\end{equation}
Condition~\eqref{eq:gumbel_tail} is the usual tail formulation of the
assumption that the distribution of $\rho$ belongs to the Gumbel
max-domain of attraction. Equivalently,
$(\max_{1\leq i\leq n}\rho_i-a_n)/b_n \ \vertk Z$,
where $Z$ has the standard Gumbel distribution.

We note that Jammalamadaka and Janson~\cite{jamsva15},
in their study of the maximum interpoint distance $D_n$, impose the
stronger condition that
\[
\mathbb{P}(\rho>a_n+b_nt)
= \frac{1+o(1)}{n}\E^{-t}
\]
holds for all $t=t_n$ satisfying $|t|\leq (d-1)\log(a_n/b_n)/2$.
As pointed out in Remark~1.2 of \cite{jamsva15},
this requirement is stronger than the ordinary Gumbel max-domain of
attraction condition~\eqref{eq:gumbel_tail}, which only requires the
above convergence for each fixed $t\in\mathbb{R}$.

Writing $\varepsilon_y$ for the Dirac measure centred at $y$, define
the point process
\[
N_n:=\sum_{i=1}^n
\varepsilon_{\left((\rho_i-a_n)/b_n,U_i\right)}
\]
on $\mathbb{R}\times\mathbb{S}^{d-1}$.

\smallskip
\noindent 
The point processes below are regarded as random elements of
$M_P(\mathbb{R}\times\mathbb{S}^{d-1})$, the space of locally finite
point measures on $\mathbb{R}\times\mathbb{S}^{d-1}$, equipped with
the vague topology.
\medskip

\begin{proposition}\label{prop:gumbel_ppp}
Suppose that \eqref{eq:gumbel_tail} holds. Then
\[
N_n \vertk N
\]
in $M_P(\mathbb{R}\times\mathbb{S}^{d-1})$, equipped with the vague
topology, where $N$ is a Poisson point process on
$\mathbb{R}\times\mathbb{S}^{d-1}$ with intensity measure
$\nu(\D s,\D u)=\E^{-s}\,\D s\,\sigma_{d-1}(\D u)$,
and $\sigma_{d-1}$ denotes the normalized surface measure on
$\mathbb{S}^{d-1}$.
\end{proposition}

\begin{proof}
Put $S_{n1}:= (\rho_1-a_n)/b_n$.
By \eqref{eq:gumbel_tail}, we have $n\,\mathbb{P}(S_{n1}>t)\rightarrow \E^{-t}$
for every $t\in\mathbb{R}$. Hence, if $s<t$, then
\[
n\,\mathbb{P}(s<S_{n1}\le t) = n\Bigl(
\mathbb{P}(S_{n1}>s)-\mathbb{P}(S_{n1}>t)
\Bigr)
\rightarrow \E^{-s}-\E^{-t}.
\]
Thus the measures $\mu_n(B):= n\,\mathbb{P}(S_{n1}\in B)$,
$B\in\mathcal{B}(\mathbb{R})$,
converge vaguely on $\mathbb{R}$ to the measure
$\mu(\D s)=\E^{-s}\,\D s$. 
Since $U_1$ is independent of $\rho_1$ and has distribution
$\sigma_{d-1}$, we have
\[
n\,\mathbb{P}\bigl((S_{n1},U_1)\in B\bigr)
=
(\mu_n\otimes\sigma_{d-1})(B),
\qquad
B\in\mathcal{B}(\mathbb{R}\times\mathbb{S}^{d-1}).
\]

Because $\sigma_{d-1}$ is a finite measure, the vague convergence
$\mu_n  \xrightarrow{v} \mu$ implies
\[
\mu_n\otimes\sigma_{d-1}
\xrightarrow{v}
\mu\otimes\sigma_{d-1}
=
\nu
\]
on $\mathbb{R}\times\mathbb{S}^{d-1}$.
Consequently,
\[
n\,\mathbb{P}\left(
\left(
\frac{\rho_1-a_n}{b_n},U_1
\right)\in\cdot
\right)
\xrightarrow{v}\nu.
\]
The standard Poisson convergence theorem for empirical point processes
therefore yields
\[
N_n=
\sum_{i=1}^n
\varepsilon_{\left((\rho_i-a_n)/b_n,U_i\right)}
\vertk N;
\]
see, for example, Proposition~3.21 of \cite{resn}.
\end{proof}

\medskip

We first note that condition~\eqref{eq:gumbel_tail} implies
\begin{equation}\label{eq:b_na_n0}
\frac{b_n}{a_n}\rightarrow 0.
\end{equation}
Indeed, let $K>0$ be fixed. By~\eqref{eq:gumbel_tail}, applied with
$t=-K$,
$n\,\mathbb{P}(\rho>a_n-Kb_n)\rightarrow \E^K$.
In particular, $\mathbb{P}(\rho>a_n-Kb_n)\rightarrow 0$.
Since $\rho\geq 0$ almost surely and, by absolute continuity,
$\mathbb{P}(\rho>0)=1$, it follows that
$a_n-Kb_n>0$ 
for all sufficiently large $n$. Hence $b_n/a_n < 1/K$ eventually. Since $K>0$ is arbitrary, \eqref{eq:b_na_n0} follows.

\medskip

Let $c \in \mathbb{R}^d$ and write $c = b_n z$ for $z \in \mathbb{R}^d$.
Moreover, let $r = a_n + b_n t$. For $X_i = \rho_i U_i$, consider the condition
$\|X_i-c\| \le r$.

The following lemma provides the key geometric approximation underlying the limit theorem. On the natural extreme-value scale, the enclosing ball problem becomes asymptotically linear: after centring at radius $a_n$ and rescaling by $b_n$, the nonlinear constraint $\|X_i-c\|\le r$
reduces to the affine inequality
\[
\frac{\rho_i-a_n}{b_n} \le t+z\cdot U_i.
\]
Thus the asymptotic behaviour of $R_n$ is determined by a functional of the limiting point process.

Here and throughout, $x\cdot y$ denotes the Euclidean inner product in $\mathbb{R}^d$.

\medskip

\begin{lemma}\label{lem:linearizationSection3}
Let $M>0$, $A>M$, and $L>0$. Then
\[
\sup_{\substack{\|z\|\le L,\,-M\le s\le A\\ u\in\mathbb{S}^{d-1}}}
\left|\frac{\|(a_n+b_ns)u-b_nz\|-a_n}{b_n}
- (s-z\cdot u) \right| \to 0 \]
as $n\to\infty$.
\end{lemma}
\medskip\noindent
Although the proof only requires $s$ to range over a compact interval,
we impose $A>M$ in order to use the same truncation parameters as in
Lemmas~\ref{lem:truncatedContinuitySection3}--\ref{lem:gumbel_truncation} and Theorem~\ref{thm:planarGumbelSection3}.

\smallskip
\begin{proof}
Write $\rho=a_n+b_ns$. Then $\|\rho u-b_nz\|^2 = (a_n+b_ns)^2
- 2b_n(a_n+b_ns)z\cdot u + b_n^2\|z\|^2$. Expanding yields
\[
\|\rho u-b_nz\|^2 =  a_n^2 + 2a_nb_n(s-z\cdot u)
+ b_n^2 \bigl( s^2-2s\,z\cdot u+\|z\|^2\bigr).
\]
Putting
\[
\delta_n  = \delta_n(s,u,z) := 2\frac{b_n}{a_n}(s-z\cdot u) + \frac{b_n^2}{a_n^2}
\bigl(s^2-2s\,z\cdot u+\|z\|^2\bigr),
\]
it follows  that $\|\rho u-b_nz\| = a_n \left(1 + \delta_n\right)^{1/2}$.
Since $s$, $u$, and $z$ range over compact sets and $b_n/a_n \to 0$,
we have
\[
\sup_{\substack{\|z\|\le L,\,-M\le s\le A\\ u\in\mathbb{S}^{d-1}}} 
|\delta_n(s,u,z)| \to 0.
\]
Hence, the expansion $(1+\delta_n)^{1/2}=1+\delta_n/2+O(\delta_n^2)$ holds uniformly over the stated ranges. 

Since the variables $s,u,z$ range over the stated compact sets,
uniformly over these sets we have $\delta_n=O\big( b_n/a_n\big)$.
Consequently,
\[
a_n\delta_n^2 =O\left(\frac{b_n^2}{a_n}\right) =o(b_n),
\]
and the second-order term in $a_n\delta_n/2$ is likewise
$O(b_n^2/a_n)=o(b_n)$. Hence
\[
\|(a_n+b_ns)u-b_nz\| = a_n+b_n(s-z\cdot u)+o(b_n)
\]
uniformly over the stated ranges. Dividing by $b_n$ yields the assertion. 
\end{proof}


\smallskip

For a locally finite point measure $m=\sum_k\varepsilon_{(s_k,u_k)}$
on $\mathbb{R}\times \mathbb{S}^{d-1}$, define 
\[
\Phi(m) :=  \inf_{z\in\mathbb{R}^d} \sup_{(s,u)\in m} (s-z\cdot u),
\]
where $\Phi(m) \in [-\infty, \infty]$. 
Note that, for the limiting Poisson point process $N$ from Proposition~\ref{prop:gumbel_ppp}
the random variable $\Phi(N)$ is finite almost surely.
Indeed, choosing $z=0$ yields
\[
\Phi(N)\le \sup_{(s,u) \in N }s ,
\]
and the right-hand side is finite almost surely, since $N((t,\infty)\times \mathbb{S}^{d-1})$ is Poisson with mean $\E^{-t}$.

On the other hand, the function
\[
z\longmapsto\sup_{(s,u)\in N}(s-z\cdot u)
\]
is coercive almost surely. To see this, fix $\delta\in(0,1/\sqrt{2})$ and choose
$v_1,\ldots,v_m\in\mathbb S^{d-1}$ such that the spherical neighbourhoods
$B(v_j,\delta) :=\{u\in\mathbb S^{d-1}:\|u-v_j\|<\delta\}$ for $j=1,\ldots,m$ 
cover $\mathbb S^{d-1}$. Put $\eta:=1-2\delta^2>0$. For each $j$, the set $\mathbb R\times B(-v_j,\delta)$  has infinite intensity and therefore contains at least one point
$(s_j,u_j)$ of $N$ almost surely. Since there are only finitely many
caps, we may choose these points simultaneously and put
$M:=\max_{1\leq j\leq m}(-s_j)_+<\infty$ almost surely.

Now let $z\neq 0$ and put $v=z/\|z\|$. Choose $j$ such that
$v\in B(v_j,\delta)$. Since $u_j\in B(-v_j,\delta)$, we have
$\|u_j+v\| \leq \|u_j+v_j\|+\|v_j-v\|<2\delta$. Consequently,
$u_j\cdot v = (\|u_j+v\|^2-2)/2 <2\delta^2-1=-\eta$,
and hence
\[
\sup_{(s,u)\in N}(s-z\cdot u) \geq s_j-z\cdot u_j \geq -M+\eta\|z\|.
\]
It follows that
\[
\sup_{(s,u)\in N}(s-z\cdot u)\longrightarrow\infty
\qquad\text{as }\|z\|\longrightarrow\infty
\]
almost surely. Furthermore,
\[
\sup_{(s,u)\in N}(s-z\cdot u) \leq \sup_{(s,u)\in N}s+\|z\|<\infty
\]
almost surely. Thus the function $z\mapsto\sup_{(s,u)\in N}(s-z\cdot u)$ 
is real-valued on $\mathbb R^d$. Being the supremum of affine
functions, it is convex and hence continuous. Together with the
preceding coercivity property, this shows that its infimum is finite
and attained. Thus $\Phi(N)\in\mathbb R$ almost surely.

\medskip
\noindent
Lemma~\ref{lem:linearizationSection3} reduces the enclosing-ball constraint to an asymptotically linear form. To apply the continuous mapping theorem to the limiting Poisson point process, it remains to establish continuity of the
associated functional under vague convergence. In what follows, the supremum over an empty set is interpreted as
$-\infty$, and extended-real-valued functionals are regarded as taking
values in $\overline{\mathbb{R}}=[-\infty,\infty]$, equipped with its
usual topology.

\medskip

\begin{lemma}\label{lem:truncatedContinuitySection3}
For $M>0$, $A>M$, and $L>0$, define
\[
\Phi_{M,A,L}(m) := \inf_{\|z\|\le L} 
\sup_{\substack{(s,u)\in m\\ -M\le s\le A}} (s-z\cdot u).
\]
Let $N$ denote the limiting Poisson point process from
Proposition~\ref{prop:gumbel_ppp}. Then, for every
$M>0$, $A>M$, and  $L>0$, the mapping $m\mapsto \Phi_{M,A,L}(m)$
is almost surely continuous at $N$ with respect to vague convergence.
\end{lemma}
\smallskip
\noindent 

\begin{proof}
Fix $L>0$, $M>0$, and $A>M$. Notice that 
$N(\{-M,A\}\times \mathbb{S}^{d-1})=0$ almost surely.
Since
\[
\nu([-M,A]\times \mathbb{S}^{d-1})
=
\int_{-M}^{A} {\rm e}^{-t}\,{\rm d}t
<
\infty,
\]
the restriction
$N_{M,A}:=N|_{[-M,A]\times \mathbb{S}^{d-1}}$
contains only finitely many points almost surely.
Define
$K_{M,A}:=
N([-M,A]\times \mathbb{S}^{d-1})$.
In what follows we work on the event of probability one on which
$K_{M,A}<\infty$
and
$N(\{-M,A\}\times \mathbb{S}^{d-1})=0$.
On this event, the number of points of
$N_{M,A}$ is finite and fixed.

Let $(m_n)_{n\ge1}$ be any sequence of locally finite point measures
on $\mathbb{R}\times\mathbb{S}^{d-1}$ such that
$m_n \to N$
vaguely as $n\to\infty$.
Note that the compact set $[-M,A]\times\mathbb{S}^{d-1}$ has boundary
$(\{-M\}\times\mathbb{S}^{d-1}) \cup \{A\}\times\mathbb{S}^{d-1})$,
which carries no mass under $N$. Hence Proposition 3.13 of Resnick \cite{resn}
implies that, for all sufficiently large $n$, $m_n([-M,A]\times \mathbb{S}^{d-1}) = K_{M,A}$. 

If $K_{M,A} =0$, then both restrictions are empty   for all sufficiently large $n$, and hence $\Phi_{M,A,L}(m_n) = \Phi_{M,A,L}(N) = - \infty$ eventually. Thus there is nothing to prove. 

Suppose now that $K_{M,A} \ge 1$. Proposition 3.13 of Resnick \cite{resn} further implies that, for all sufficiently large $n$, 
the points of $m_n$ in $[-M,A]\times\mathbb{S}^{d-1}$ 
can be labelled as
\[
(s_1^{(n)},u_1^{(n)}),\ldots,
(s_{K_{M,A}}^{(n)},u_{K_{M,A}}^{(n)}),
\]
in such a way that $(s_j^{(n)},u_j^{(n)}) \to (s_j,u_j)$
for $j=1,\ldots,K_{M,A}$. 
Define, for $\|z\|\le L$,
\[
H(z) := \max_{1\le j\le K_{M,A}} (s_j-z\cdot u_j), \quad \text{and} \quad H_n(z)
:= \max_{1\le j\le K_{M,A}} (s_j^{(n)}-z\cdot u_j^{(n)}).
\]
Then $\sup_{\|z\|\le L}|H_n(z)-H(z)|\to 0$. Indeed, for $\|z\|\le L$,
\[
|H_n(z)-H(z)| \le \max_{1\le j\le K_{M,A}} \left( |s_j^{(n)}-s_j| +  L\|u_j^{(n)}-u_j\| \right),
\]
and the right-hand side converges to $0$ as $n\to\infty$. Consequently,
\[
\left|\inf_{\|z\|\le L}H_n(z) - \inf_{\|z\|\le L}H(z) \right|
\le \sup_{\|z\|\le L}|H_n(z)-H(z)| \to 0.
\]
Since $\Phi_{M,A,L}(m_n) = \inf_{\|z\|\le L}H_n(z)$ and
$\Phi_{M,A,L}(N) = \inf_{\|z\|\le L}H(z)$,
it follows that $\Phi_{M,A,L}(m_n) \to \Phi_{M,A,L}(N)$, which proves the claim.
\end{proof}

To complete the proof of Theorem~\ref{thm:planarGumbelSection3}, it remains to control the optimization over the centre variable on the natural scale $b_n$.
The next lemma shows that the centres of the exact enclosing balls are
tight on this scale and that, with arbitrarily large probability, the
linearized optimization problem may be restricted to a sufficiently
large compact subset of $\mathbb{R}^d$.

\begin{lemma}\label{lem:tight-centres}
Let $C_n$ denote the centre of the smallest enclosing ball of
$X_1,\ldots,X_n$, and put
\[ 
S_{ni}:=\frac{\rho_i-a_n}{b_n}, \qquad i=1,\ldots,n,
\]
where $\rho_i = \|X_i\|$. 
Assume that \eqref{eq:gumbel_tail} holds. Then $C_n/b_n=O_\mathbb{P}(1)$.
Moreover, define 
\begin{equation}\label{defgnz}
G_n(z):=\max_{1\le i\le n}(S_{ni}-z\cdot U_i), \qquad z\in\mathbb{R}^d.
\end{equation}
Then, for every $\varepsilon>0$, there exists $L>0$ such that, for
all sufficiently large $n$,
\[
\mathbb P\left(
\inf_{z\in\mathbb R^d}G_n(z)
=
\inf_{\|z\|\leq L}G_n(z)
\right)
\geq 1-\varepsilon.
\]

\end{lemma}

\smallskip
\noindent 

\begin{proof}
Fix $\varepsilon>0$. We first establish the compactness property for
the linearized functional $G_n$. Choose $\delta\in(0,1/\sqrt{2})$ and
let $v_1,\ldots,v_m\in\mathbb S^{d-1}$ be such that the spherical neighbourhoods
\[
B(v_j,\delta):=\{u\in\mathbb S^{d-1}:\|u-v_j\|<\delta\},
\qquad j=1,\ldots,m,
\]
cover $\mathbb S^{d-1}$. For $M>0$ and $j=1,\ldots,m$, let
$A_{n,j}:= \left\{ N_n\bigl([-M,\infty)\times B(-v_j,\delta)\bigr)\geq 1 \right\}$.
Thus, on $A_{n,j}$, there exists an index $i\leq n$ such that
$S_{ni}\geq -M$ and $U_i\in B(-v_j,\delta)$.

For fixed $j$, the random variable
$N_n\bigl([-M,\infty)\times B(-v_j,\delta)\bigr)$
has a binomial distribution with parameters $n$ and
$p_{n,j} = \mathbb P(S_{n1}\geq -M)\,
\sigma_{d-1}\bigl(B(-v_j,\delta)\bigr)$.
Since
$n\mathbb P(S_{n1}\geq -M)\rightarrow {\rm e}^M$,
we have
$np_{n,j}\rightarrow {\rm e}^M\sigma_{d-1}\bigl(B(-v_j,\delta)\bigr)$,
and consequently
$\mathbb P(A_{n,j}^c) \rightarrow \exp\left\{-\E^M\sigma_{d-1}\bigl(B(-v_j,\delta)\bigr)
\right\}$.
Since there are only finitely many neighbourhoods, the union bound shows that
$M$ can be chosen so large that, with
$E_n:=\bigcap_{j=1}^m A_{n,j}$, we have
\[
\liminf_{n\to\infty}\mathbb P(E_n)\geq 1-\frac{\varepsilon}{4}
\]
and thus $\mathbb P(E_n)\geq 1-\varepsilon/2$ for sufficiently large $n$.

Put $\eta:=1-2\delta^2>0$. On the event $E_n$, the following property holds. For every
$v\in\mathbb S^{d-1}$, choose $j$ such that $v\in B(v_j,\delta)$. By the definition of $E_n$, there exists an index $i=i(v)$ such that $S_{ni}\geq -M$ and $U_i\in B(-v_j,\delta)$. Hence  $\|U_i+v\| \leq \|U_i+v_j\|+\|v_j-v\| <2\delta$.
Since $\|U_i\|=\|v\|=1$, we have $2+2U_i\cdot v=\|U_i+v\|^2<4\delta^2$ and therefore
$U_i\cdot v\leq -\eta$.

Now let $z\in\mathbb R^d$, $z\neq0$, and write
$v=z/\|z\|$. Choosing $i=i(v)$ as above gives
\[
S_{ni}-z\cdot U_i = S_{ni}-\|z\|(U_i\cdot v) \geq -M+\eta\|z\|.
\]
Consequently, on $E_n$, $G_n(z)\geq -M+\eta\|z\|$ for each
$z\neq 0$. On the other hand,
\[
G_n(0)=\max_{1\leq i\leq n}S_{ni}=O_{\mathbb P}(1).
\]
Hence we may choose $K>0$ such that $\mathbb P\bigl(G_n(0)>K\bigr)\leq \varepsilon/2$
for all sufficiently large $n$.

Choose $L>(M+K)/\eta$. Then, on the event $E_n\cap\{G_n(0)\leq K\}$,
we have, for every $\|z\|>L$, 
\[
G_n(z)\geq -M+\eta\|z\|>K\geq G_n(0).
\]
Since $0\in\{z:\|z\|\leq L\}$, it follows that
$\inf_{z\in\mathbb R^d}G_n(z) = \inf_{\|z\|\leq L}G_n(z)$
on this event. Therefore, for all sufficiently large $n$,
\[
\mathbb P\left(\inf_{z\in\mathbb R^d}G_n(z) =
\inf_{\|z\|\leq L}G_n(z) \right) \geq 1-\varepsilon.
\]

We now turn to the exact enclosing-ball centres. Define
\[
H_n(z):= \frac{1}{b_n} \left(\max_{1\leq i\leq n}\|\rho_iU_i-b_nz\|-a_n \right),
\qquad z\in\mathbb R^d.
\]
Then
\[
\frac{R_n-a_n}{b_n} = \inf_{z\in\mathbb R^d}H_n(z),
\]
and the minimizer of $H_n$ is $C_n/b_n$.

We use the same event $E_n$ as above. Thus, on $E_n$, for every
$v\in\mathbb S^{d-1}$ there exists an index $i=i(v)$ such that
$S_{ni}\geq -M$ and $U_i\cdot v\leq-\eta$.
Let $z\neq0$ and write $v=z/\|z\|$. Since
$\rho_i=a_n+b_nS_{ni}$, we have
\begin{align*}
\|\rho_iU_i-b_nz\|^2
&= \rho_i^2 -2b_n\|z\|\rho_i(U_i\cdot v) +b_n^2\|z\|^2 \\
&\geq \rho_i^2 +2\eta b_n\|z\|\rho_i +b_n^2\|z\|^2 \\
&\geq \bigl(\rho_i+\eta b_n\|z\|\bigr)^2, 
\end{align*}
where the last inequality follows from $0<\eta\leq1$.
Hence
\[
\|\rho_iU_i-b_nz\| \geq \rho_i+\eta b_n\|z\|
= a_n+b_n\bigl(S_{ni}+\eta\|z\|\bigr),
\]
and therefore $H_n(z)\geq -M+\eta\|z\|$ for all $z\neq0$ on $E_n$.
Moreover,
\[
H_n(0)
=
\frac{1}{b_n}
\left(
\max_{1\leq i\leq n}\rho_i-a_n
\right)
=
\max_{1\leq i\leq n}S_{ni}
=
G_n(0).
\]
Thus, on $E_n\cap\{H_n(0)\leq K\}$, we have
$H_n(z)\geq -M+\eta\|z\|>K\geq H_n(0)$ for every $\|z\|>L$.
Consequently, the minimizer $C_n/b_n$ of $H_n$ satisfies
$\|C_n\|/b_n\leq L$ on this event. Since
\[
\liminf_{n\to\infty}\mathbb P(E_n)\geq1-\frac{\varepsilon}{2}
\]
and $H_n(0)=O_{\mathbb P}(1)$, it follows that
$C_n/b_n=O_{\mathbb P}(1)$, which completes the proof.
\end{proof}

We are now in a position to derive the asymptotic representation of the
minimal enclosing radius. Proposition~\ref{prop:gumbel_ppp}
describes the asymptotic extremal configuration through a Poisson point
process, while Lemma~\ref{lem:linearizationSection3} shows that the
geometric covering constraint becomes asymptotically linear on the
extreme-value scale. The limiting fluctuations of $R_n$ are therefore
governed by the variational functional $\Phi$.

We next show that the truncations in the extreme-value coordinate
$s$ introduced in Lemma~\ref{lem:truncatedContinuitySection3} are asymptotically negligible. More precisely, for fixed compact sets in the optimization variable
$z$, points with very large or very small values of $S_{ni}$
do not contribute to the variational functional with high
probability. To this end, recall the functional $G_n(z)$ 
introduced in \eqref{defgnz}. For $M>0$ and $A>M$, define the truncated functional
\[
G_{n,M,A}(z) = \max_{\substack{1\le i\le n\\ -M\le S_{ni}\le A}}
(S_{ni}-z\cdot U_i). \]
Similarly, for the limiting Poisson process $N$, define
\[
G(z) = \sup_{(s,u)\in N}(s-z\cdot u), \qquad  G_{M,A}(z)
= \sup_{\substack{(s,u)\in N\\ -M\le s\le A}} (s-z\cdot u).
\]
In what follows, maxima over empty index sets are interpreted as $-\infty$.

\begin{lemma}\label{lem:gumbel_truncation}
Fix $L>0$. Then
\[
\lim_{M\to\infty}\liminf_{A\to\infty}\liminf_{n\to\infty}
\mathbb{P}\Bigl(G_n(z)=G_{n,M,A}(z) \text{ for all }\|z\|\le L \Bigr)=1.
\]
Moreover,
\[
\lim_{M\to\infty}\liminf_{A\to\infty}
\mathbb{P}\Bigl(G(z)=G_{M,A}(z) \text{ for all }\|z\|\le L \Bigr)=1.
\]
\end{lemma}

\begin{proof}
Fix $L>0$, let $M>4L$, and let $A>M$. Put
\[
U_{n,A}:=\left\{\max_{1\le i\le n}S_{ni}>A\right\}, \qquad 
B_{n,M}:= \left\{\max_{1\le i\le n}S_{ni}\ge -\frac{M}{2}\right\}.
\]
On $U_{n,A}^{c}$, no point with $S_{ni}>A$ occurs, so the upper
truncation has no effect.

On $B_{n,M}$ there exists an index $i$ such that
$S_{ni}\ge -M/2$. On $U_{n,A}^{c}\cap B_{n,M}$ this index also
satisfies $S_{ni}\le A$, and hence it is retained in
$G_{n,M,A}$. For every $\|z\|\le L$,
$S_{ni}-z\cdot U_i\ge - M/2-L$.
On the other hand, if $S_{nj}<-M$, then
$S_{nj}-z\cdot U_j<-M+L$.
Since $M>4L$, we have
$-M/2-L>-M+L$.
Thus no point with $S_{nj}<-M$ can attain the maximum, uniformly
for $\|z\|\le L$. Consequently, on
$U_{n,A}^{c}\cap B_{n,M}$, we have
$G_n(z)=G_{n,M,A}(z)$ for all $\|z\|\le L$.

By \eqref{eq:gumbel_tail}, $\mathbb{P}(U_{n,A}) \rightarrow 1-\exp(-\E^{-A})$ and
\[
\mathbb{P}(B_{n,M}^{c})
= \mathbb{P}\left(
\max_{1\le i\le n}S_{ni}<-\frac{M}{2}
\right) \rightarrow \exp(-\E^{M/2}).
\]
Hence
\[
\liminf_{n\to\infty}
\mathbb{P}\Bigl(G_n(z)=G_{n,M,A}(z) \text{ for all }\|z\|\le L
\Bigr) \ge
1-\bigl(1-\E^{-\E^{-A}}\bigr)-\E^{-\E^{M/2}}.
\]
Letting first $A\to\infty$ and then $M\to\infty$ proves the first
assertion.

For the limiting Poisson process, put
\[
U_A:=\{N((A,\infty)\times\mathbb S^{d-1})>0\}, \qquad 
B_M:=
\{N([-M/2,\infty)\times\mathbb S^{d-1})\ge1\}.
\]
The same deterministic comparison shows that, on
$U_A^c\cap B_M$, we have $G(z)=G_{M,A}(z)$ for all $\|z\|\le L$.  Moreover,
$\mathbb{P}(U_A)=1-\exp(-\E^{-A})$ and
$\mathbb{P}(B_M^c)=\exp(-\E^{M/2})$.
Letting first $A\to\infty$ and then $M\to\infty$ proves the second
assertion.
\end{proof}

\medskip

We are now in a position to derive the full Poisson functional
limit theorem for the enclosing radius. Proposition~\ref{prop:gumbel_ppp} provides the convergence of the extremal point process, Lemma~\ref{lem:linearizationSection3} yields
the asymptotic linearization of the geometric constraint, while
Lemmas~\ref{lem:truncatedContinuitySection3}--\ref{lem:gumbel_truncation} establish continuity, coercivity, and the asymptotic negligibility of the truncations. Combining these
ingredients yields the following limit theorem.

\begin{theorem}\label{thm:planarGumbelSection3}
Under \eqref{eq:gumbel_tail},
\[
\frac{R_n-a_n}{b_n} \vertk \Phi(N)
= \inf_{z\in\mathbb{R}^d} \sup_{(s,u)\in N} (s-z\cdot u),
\]
where $N$ is the Poisson point process from Proposition~\ref{prop:gumbel_ppp}
with intensity measure $\nu(\D s,\D u)=\E^{-s}\,\D s\,\sigma_{d-1}(\D u)$.
\end{theorem}

\smallskip

\begin{proof} Put $S_{ni}:=(\rho_i-a_n)/b_n$ for $i=1,\ldots,n$, and define
\[
H_n(z):= \frac{1}{b_n} \left(\max_{1\le i\le n}\|\rho_iU_i-b_nz\|-a_n\right), \qquad z \in \mathbb{R}^d,
\]
so that
\[
\frac{R_n-a_n}{b_n} = \inf_{z\in\mathbb{R}^d}H_n(z).
\]
Furthermore, let
\[
G_n(z):=\max_{1\le i\le n}(S_{ni}-z\cdot U_i).
\]

We first restrict the optimization to compact sets in the $z$-variable.
By Lemma~\ref{lem:tight-centres}, both the minimization of $H_n$ and $G_n$ may, with arbitrarily large probability, be restricted to a sufficiently large ball. Hence, for every $\varepsilon>0$ there exists $L>0$ such that, for all sufficiently large $n$,
\[
\mathbb{P}\left( \inf_{z\in\mathbb{R}^d}H_n(z)
= \inf_{\|z\|\le L}H_n(z)\right)\ge 1-\varepsilon
\]
and
\[
\mathbb{P}\left(\inf_{z\in\mathbb{R}^d}G_n(z)
= \inf_{\|z\|\le L}G_n(z)\right) \ge 1-\varepsilon.
\]
It is therefore enough to prove convergence of the compactly restricted functionals and then let $L\to\infty$. To this end, define
\[
G_{n,M,A,L} := \inf_{\|z\|\le L}
\max_{\substack{1\le i\le n\\ -M\le S_{ni}\le A}} (S_{ni}-z\cdot U_i),
\]
where $M>0$, $A>M$, and $L>0$. Note that $G_{n,M,A,L}=\Phi_{M,A,L}(N_n)$.
By Proposition~\ref{prop:gumbel_ppp} and Lemma~\ref{lem:truncatedContinuitySection3}, 
for all $M>0$ and $A>M$,
\[
G_{n,M,A,L}  = \Phi_{M,A,L}(N_n) \vertk \Phi_{M,A,L}(N).
\]
Moreover, Lemma~\ref{lem:gumbel_truncation} implies
\[
\lim_{M\to\infty}\liminf_{A\to\infty}\liminf_{n\to\infty}
\mathbb P\left(
\inf_{\|z\|\le L}G_n(z)=G_{n,M,A,L}
\right)=1.
\]
The same lemma gives
\[
\lim_{M\to\infty}\liminf_{A\to\infty}
\mathbb P\left(
\Phi_L(N)=\Phi_{M,A,L}(N)
\right)=1,
\]
where
\[
\Phi_L(N):=\inf_{\|z\|\le L}\sup_{(s,u)\in N}(s-z\cdot u).
\]
Choosing sequences $M_k\to\infty$ and $A_k\to\infty$ with
$A_k>M_k$ such that the two approximation probabilities above tend
to one, we may apply the converging-together theorem 
(see, e.g., Theorem 3.2 of \cite{bi})
with
$G_{n,M_k,A_k,L}$ as the $k$th approximation.
Consequently, 
\[
\inf_{\|z\|\le L}G_n(z) \vertk \Phi_L(N).
\]

We now compare the exact and the linearized compact functionals.
Put
\[
h_{ni}(z):= \frac{\|\rho_iU_i-b_nz\|-a_n}{b_n}, \qquad 1\leq i\leq n,
\]
and define
\[
H_{n,M,A}(z):= \max_{\substack{1\leq i\leq n\\ -M\leq S_{ni}\leq A}} h_{ni}(z),
\qquad \|z\|\leq L,
\]
where, as before, a maximum over an empty index set is interpreted as
$-\infty$.

For $\|z\|\leq L$, the triangle inequality gives
$\rho_i-b_nL \leq \|\rho_iU_i-b_nz\| \leq \rho_i+b_nL$,
and hence
$S_{ni}-L\leq h_{ni}(z)\leq S_{ni}+L$.
For $M>4L$ and $A>M$, put
\[
U_{n,A}:=\left\{\max_{1\leq i\leq n}S_{ni}>A\right\},
\qquad
B_{n,M}:=\left\{\max_{1\leq i\leq n}S_{ni}\geq-\frac M2\right\}.
\]
On $U_{n,A}^c$, no point with $S_{ni}>A$ occurs, so the upper
truncation has no effect. On $B_{n,M}$ there exists an index $i$
such that $S_{ni}\geq-M/2$. Thus, on
$U_{n,A}^c\cap B_{n,M}$, this index satisfies
$-M/2\leq S_{ni}\leq A$, and for every $\|z\|\leq L$, we have
$h_{ni}(z)\geq- M/2-L$. On the other hand, if $S_{nj}<-M$, then
$h_{nj}(z)<-M+L$, if $\|z\|\leq L$.
Since $M>4L$, we have
$-M/2-L>-M+L$. Consequently, on $U_{n,A}^c\cap B_{n,M}$,
\[
H_n(z)=H_{n,M,A}(z) \quad\text{and}\quad G_n(z)=G_{n,M,A}(z)
\]
simultaneously for all $\|z\|\leq L$.
For fixed $M,A,L$, define
\[
\Delta_n(M,A,L):= \sup_{\substack{\|z\|\leq L,\,-M\leq s\leq A\\
u\in\mathbb S^{d-1}}} \left| \frac{\|(a_n+b_ns)u-b_nz\|-a_n}{b_n}
-(s-z\cdot u) \right|.
\]
By Lemma~\ref{lem:linearizationSection3}, we have
$\Delta_n(M,A,L)\rightarrow 0$.
On $U_{n,A}^c\cap B_{n,M}$ the truncated index set is non-empty, and
therefore
\[
\left| \inf_{\|z\|\leq L}H_n(z) - \inf_{\|z\|\leq L}G_n(z) \right|
\leq \Delta_n(M,A,L).
\]
Moreover, by~(3.1), we have
$\mathbb P(U_{n,A})\rightarrow 1-\exp(-\E^{-A})$
and
\[
\mathbb P(B_{n,M}^c) = \mathbb P\left(\max_{1\leq i\leq n}S_{ni}<-\frac M2\right) \rightarrow \exp(-e^{M/2}).
\]
It follows that, for every $\varepsilon>0$,
\[
\limsup_{n\to\infty}
\mathbb P\left(\left| \inf_{\|z\|\leq L}H_n(z) -
\inf_{\|z\|\leq L}G_n(z) \right|>\varepsilon \right)
\leq 1-\exp(-\E^{-A})+\exp(-\E^{M/2}).
\]
Letting first $A\to\infty$ and then $M\to\infty$ gives
\[
\inf_{\|z\|\leq L}H_n(z) - \inf_{\|z\|\leq L}G_n(z) \stk 0.
\]
Together with
$\inf_{\|z\|\leq L}G_n(z) \vertk \Phi_L(N)$
this yields
$\inf_{\|z\|\leq L}H_n(z) \vertk \Phi_L(N)$ 
by Slutsky's lemma. Finally let $L\to\infty$. Since
\[
\Phi_L(N) = \inf_{\|z\|\leq L} \sup_{(s,u)\in N}(s-z\cdot u)
\downarrow \inf_{z\in\mathbb R^d} \sup_{(s,u)\in N}(s-z\cdot u) = \Phi(N)
\]
almost surely, we have $\Phi_L(N)\vertk \Phi(N)$.
Moreover, by Lemma~\ref{lem:tight-centres}, for every $\varepsilon>0$ there exists
$L>0$ such that, for all sufficiently large $n$,
\[
\mathbb P\left(\inf_{z\in\mathbb R^d}H_n(z) = \inf_{\|z\|\leq L}H_n(z) \right)
\geq 1-\varepsilon.
\]
Hence another application of the converging-together theorem yields
\[
\frac{R_n-a_n}{b_n} = \inf_{z\in\mathbb R^d}H_n(z) \vertk \Phi(N).
\]
This proves the theorem.
\end{proof}

\medskip

%

\begin{remark}\label{rem:gumbel-planar}
%
In dimension $d=2$, the distribution of the limiting functional
$\Phi(N)$ can be identified explicitly. Indeed, in the planar Gaussian
case, Theorems~\ref{haupt} and~\ref{thm:planarGumbelSection3}, together with
\[
a_n=\sqrt{2\log n},
\qquad
b_n=(2\log n)^{-1/2},
\]
yield
\[
\mathbb P\bigl(\Phi(N)\leq t\bigr)
=
G(t)\bigl(1+\E^{-t}\bigr),
\qquad t\in\mathbb R,
\]
where $G(t)=\exp(-\E^{-t})$ is the Gumbel distribution function.

Since the limiting Poisson point process in Theorem~\ref{thm:planarGumbelSection3} has intensity
measure $\E^{-s}\,\D s\,\sigma_1(\D u)$ for every radial distribution satisfying~\eqref{eq:gumbel_tail}, the distribution of $\Phi(N)$ is universal within the planar Gumbel regime. Consequently, for every planar spherically symmetric distribution satisfying~\eqref{eq:gumbel_tail}, with corresponding norming sequences $a_n$ and $b_n$,
\[
\lim_{n\to\infty}
\mathbb P\left(\frac{R_n-a_n}{b_n}\leq t\right)
=
G(t)\bigl(1+\E^{-t}\bigr),
\qquad t\in\mathbb R.
\]
Thus, in dimension two, the Poisson functional limit of
Theorem~\ref{thm:planarGumbelSection3} coincides with the limit law arising from Daniels' classical Gaussian formula. In dimensions $d\geq3$, an explicit
characterization of the distribution of $\Phi(N)$ appears to remain
open.
\end{remark}

\medskip

%
\begin{remark}\label{sec3:rem38}
%
As an immediate consequence of Theorem~\ref{thm:planarGumbelSection3},
we obtain a limit theorem for the Kotz family of spherically symmetric
distributions (see, e.g., Section 3.2 of \cite{fkn}). Suppose that $X$ has density
\[
f(x)
=
K \|x\|^{2(q-1)}
\exp\{-r\|x\|^{2s}\},
\qquad x\in\mathbb R^d,
\]
where $r>0$, $s>0$, $q>1-d/2$, and $K>0$ is the normalizing
constant. 

Writing $\rho=\|X\|$ and denoting by $|\mathbb S^{d-1}|$ the
surface area of the unit sphere, polar coordinates give
\[
\mathbb P(\rho>t)
=
K|\mathbb S^{d-1}|
\int_t^\infty
y^{2q+d-3}\exp\{-r y^{2s}\}\,dy .
\]
With the substitution $v=ry^{2s}$, this becomes
\[
\mathbb P(\rho>t)
= \frac{K|\mathbb S^{d-1}|}{2s}\,
r^{-(2q+d-2)/(2s)} \Gamma\left(\frac{2q+d-2}{2s},\,rt^{2s}
\right),
\]
where $\Gamma(a,t)$ denotes the upper incomplete gamma function.
Using  $\Gamma(a,t)\sim t^{a-1} \E^{-t}$ as $t\to\infty$,
we obtain $\mathbb P(\rho>t) \sim C\,t^\beta\exp\{-rt^{2s}\}$ as $t\to\infty$,
where $\beta=2q+d-2-2s$ and $C=K|\mathbb S^{d-1}|/(2sr)>0$.

Let $a_n\to\infty$ be chosen such that
\[
nC a_n^\beta \exp\{-r a_n^{2s}\}\to1,
\]
and put
\[
b_n=\frac{1}{2sr\,a_n^{2s-1}}.
\]
Then $b_n/a_n\to0$, and, for every fixed $t\in\mathbb R$,
\[
\frac{\mathbb P(\rho>a_n+b_nt)}
     {\mathbb P(\rho>a_n)}
\rightarrow {\rm e}^{-t}.
\]
Indeed, $(1+b_nt/a_n)^\beta\to1$, while
$r((a_n+b_nt)^{2s}-a_n^{2s})
 =2sr a_n^{2s-1}b_nt+o(1)\to t$.
Consequently, $n\mathbb P(\rho>a_n+b_nt)\to \E^{-t}$ for every
$t\in\mathbb R$, so that condition~\eqref{eq:gumbel_tail} holds.
Hence Theorem~\ref{thm:planarGumbelSection3}
gives 
\[ \frac{R_n-a_n}{b_n} \vertk \Phi(N).
\]
More precisely, the centring constants may be chosen such that
\[
a_n^{2s}
= \frac{1}{r} \left( \log n + \frac{\beta}{2s}\log\log n + \log C
- \frac{\beta}{2s}\log r + o(1) \right),
\]
and the corresponding scale is $b_n=(2sr a_n^{2s-1})^{-1}$.

The centred Gaussian distribution with identity covariance matrix is
contained as the special case $q=1$, $s=1$, and $r=1/2$. In this case
$\beta=d-2$,
\[
a_n = \sqrt{2\log n} + \frac{(d-2)\log\log n}{2\sqrt{2\log n}}
- \frac{\log\Gamma(d/2)}{\sqrt{2\log n}} + o((\log n)^{-1/2}),
\]
and $b_n\sim(2\log n)^{-1/2}$. Thus Theorem~\ref{thm:planarGumbelSection3} includes, in particular,
the Gaussian case in arbitrary dimension.
\end{remark}

%
%
\begin{remark}
It is instructive to compare the normalization in
Theorem~\ref{thm:planarGumbelSection3} with the corresponding
normalization for the diameter $D_n$ of a symmetric Kotz sample
obtained by Henze and Klein~\cite{hekl96}.
Consider the special case of Remark~\ref{sec3:rem38} with
$q=b$, $s=1$, and $r=\eta$. Thus
\[
f(x)\propto \|x\|^{2(b-1)}  \exp(-\eta\|x\|^2),
\]
and the radial tail satisfies
\[
\mathbb P(\|X\|>t) \sim C t^{d+2b-4}\exp(-\eta t^2),
\qquad t\to\infty,
\]
for some constant $C>0$. Hence the centring constants $a_n$
in Theorem~\ref{thm:planarGumbelSection3} may be chosen such that
\[
a_n^2 = \frac{1}{\eta} \left( \log n+ \frac{1}{2}(d+2b-4)\log\log n+O(1)
\right),
\]
and the corresponding scale satisfies
\[
b_n=\frac{1}{2\eta a_n} \sim \frac{1}{2\sqrt{\eta\log n}}.
\]
By contrast, the normalization for $D_n/2$ in~\cite{hekl96} is based
on centring constants $r_1=r_1(n)$ satisfying
\[
r_1^2 = \frac{1}{\eta}
\left(\log n+ \frac{1}{4}(d+4b-7)\log\log n +O(\log\log\log n) \right).
\]
Consequently,
\[
a_n-r_1 \sim \frac{d-1}{8\sqrt{\eta\log n}}\log\log n, \qquad n\to\infty.
\]
In particular, $a_n-r_1\to0$, whereas
$(a_n-r_1)/b_n\sim (d-1)\log\log n/4\to\infty$.
Thus, although the two centring constants are asymptotically
indistinguishable on the original scale, their difference diverges
when measured on the fluctuation scale of $R_n$.
\end{remark}


\section{The Fréchet case}\label{secfrechet}


In this section we consider the heavy–tailed (Fréchet) regime for the radial distribution.
We assume that this distribution has a regularly varying tail with index
$-\alpha$, that is
\begin{equation}\label{eq-reg-var}
1-F(s)=s^{-\alpha}L(s), \qquad s\to\infty,
\end{equation}
where $\alpha>0$, $L$ is slowly varying, and $F$ denotes the
distribution function of $\|X\|$ (see, e.g., Appendix A3 of \cite{ekm97}).

In this setting the asymptotic behaviour of extreme observations is
conveniently described by point–process methods. Roughly speaking,
after suitable normalization the extreme sample points converge to a
Poisson point process. Geometric functionals such as the diameter or
the enclosing radius can then be analysed via the geometry of this
limiting Poisson cloud.

To make this precise we adopt the notation of Chapter~3 of
Resnick~\cite{resn}. 
Let $E=[-\infty,\infty]^d\setminus\{0\}$
denote the punctured compactification of $\mathbb R^d$, equipped with
the relative topology and its Borel $\sigma$-field $\mathcal E$.
 Compact
subsets of $E$ are those closed sets that are bounded away from the
origin.
 Let
$M_P(E)$ denote the set of point measures of the form $m= \sum_{j=1}^\infty \varepsilon_{x_j}$,
where $\{x_j: j \ge 1\}$ is a countable collection of points of $E$, and $m$ takes finite values
on compact subsets of $E$. 
The space $M_P(E)$ is equipped with the vague topology and its
Borel $\sigma$-field.

A point process on $E$ is a random element of $M_P(E)$ that is defined on some suitable probability space. A Poisson point process on $E$ with intensity measure $\nu$ is a point
process $N=\sum_{k\ge1}\varepsilon_{P_k}$ such that for each Borel set
$B\subset E$ with $\nu(B) < \infty$, the random variable $N(B)$ has a Poisson distribution with mean $\nu(B)$, and the counts on disjoint sets are independent. Whenever a point measure is used in a set-theoretic or geometric
operation, we identify it with its support.

For a bounded (possibly infinite) subset $A$ of $\mathbb R^d$, let
\[
\mathcal C(A)=\inf_{y\in\mathbb R^d}\sup_{x\in A}\|x-y\|
\]
denote the circumradius (radius of the smallest enclosing ball) of
$A$. Thus, $R_n=\mathcal C(\{X_1,\ldots,X_n\})$. Notice that
$\mathcal C(A)=\mathcal C(\overline{A})$, since, for every
$y\in\mathbb R^d$, continuity of $x\mapsto\|x-y\|$ implies
$\sup_{x\in A}\|x-y\|=\sup_{x\in\overline{A}}\|x-y\|$.

To state the main theorem, let
\[
c_n := \inf\{t >0: F(t) \ge 1-  1/n\}.
\]
Then $c_n\to\infty$ and $n\mathbb{P}(\|X\|>c_n)\to 1$. Define
\[
N_n=\sum_{k=1}^n \varepsilon_{c_n^{-1}X_k},
\]
and let $N=\sum_{k\ge1}\varepsilon_{P_k}$ be a Poisson point process  
 with intensity measure
\begin{equation}\label{sec-4-intmeas}
\nu\left(\bigg{\{}x: \|x\| >r, \frac{x}{\|x\|} \in A\bigg{\}}\right) = r^{-\alpha} \sigma_{d-1}(A), \quad r >0, A \in {\cal B}(\mathbb{S}^{d-1}),
\end{equation}
where, as before,  $\sigma_{d-1}$ denotes the uniform probability measure on the $\sigma$-field ${\cal B}(\mathbb{S}^{d-1})$ of Borel sets of $\mathbb{S}^{d-1}$.

\smallskip
\noindent
Notice that the limiting Poisson cloud is bounded almost surely. Indeed,
for every $r>0$,
\[
\mathbb{P}\left(\sup_{k\geq1}\|P_k\|\leq r\right)
=
\mathbb{P}\left(N(\{x:\|x\|>r\})=0\right)
=
\exp(-r^{-\alpha}),
\]
which tends to $1$ as $r\to\infty$. In particular,
$C(\{P_k:k\geq1\})<\infty$ almost surely.
A technical difficulty in the Fr\'{e}chet setting is that the limiting Poisson cloud contains infinitely many points accumulating near the origin. To obtain continuity of the
convex-hull and circumradius functionals under vague convergence, we therefore first truncate the point processes away from the origin and subsequently show that the omitted points are
asymptotically negligible.

\begin{lemma}\label{lem:Frechet-truncation}
For $\varepsilon>0$, define the truncated point processes
\[
N_n^{(\varepsilon)} := \sum_{k=1}^n \varepsilon_{c_n^{-1}X_k}
1\{\|X_k\|>\varepsilon c_n\} \quad \text{and} \quad  N^{(\varepsilon)} :=
\sum_{k\ge1} \varepsilon_{P_k} 1\{\|P_k\|>\varepsilon\}.
\]
Then, for every fixed $\varepsilon>0$,
$N_n^{(\varepsilon)} \vertk N^{(\varepsilon)}$ as $n \to \infty$ in $M_P(E)$, endowed with the vague topology. 

Moreover, the corresponding circumradii satisfy
\[
\mathcal{C}(\operatorname{conv}(\{0\}\cup N_n^{(\varepsilon)}))
\vertk \mathcal{C}(\operatorname{conv}(\{0\}\cup N^{(\varepsilon)})) \ \text{as} \  n\to\infty.
\]
Finally, for every $\delta>0$,
\[
\lim_{\varepsilon\downarrow0} \limsup_{n\to\infty} \mathbb{P}\!\left( \left| c_n^{-1}R_n
- \mathcal{C}(\operatorname{conv}(\{0\}\cup N_n^{(\varepsilon)}))\right| >\delta \right) =0.
\]
The analogous statement also holds for the limiting Poisson cloud:
\[
\lim_{\varepsilon\downarrow 0} \mathbb{P}\!\left( \left|
\mathcal{C}(\{0\}\cup\{P_k:k\ge1\}) -
\mathcal{C}(\operatorname{conv}(\{0\} \cup N^{(\varepsilon)}))
\right| >\delta \right) =0. \]
\end{lemma}

\smallskip
\noindent 
\begin{proof}
By regular variation and the definition of $c_n$, for every $r>0$
and every Borel set $A\subset \mathbb{S}^{d-1}$ satisfying
$\sigma_{d-1}(\partial A)=0$, we have
\[
\begin{aligned}
n\mathbb{P}\left( \|c_n^{-1}X\|>r,\,  \frac{X}{\|X\|}\in A \right)
& = n\mathbb{P}(\|X\|>rc_n)\sigma_{d-1}(A)\\
& = n\mathbb{P}(\|X\|>c_n) \frac{\mathbb{P}(\|X\|>rc_n)}      {\mathbb{P}(\|X\|>c_n)}
\sigma_{d-1}(A)\\
&\rightarrow r^{-\alpha}\sigma_{d-1}(A).
\end{aligned}
\]
Consequently,
$n\mathbb{P}(c_n^{-1}X\in\cdot) \xrightarrow{v}\nu$ on $E$, where $\nu$ is given by \eqref{sec-4-intmeas}. Hence, by Proposition~3.21 of Resnick~\cite{resn},
$N_n\vertk N$ in $M_P(E)$ with respect to the vague topology.

Fix $\varepsilon>0$, and put $E_\varepsilon:=\{x\in E:\|x\|>\varepsilon\}$.
Since $N\bigl(\{x:\|x\|=\varepsilon\}\bigr)=0$ almost surely,
the restriction mapping $m\mapsto m|_{E_\varepsilon}$ is continuous
at $N$ almost surely. Hence $N_n^{(\varepsilon)} \vertk N^{(\varepsilon)}$.

Moreover, $N^{(\varepsilon)}$ contains only finitely many points
almost surely, since
$\nu\bigl(\{x:\|x\|>\varepsilon\}\bigr)=\varepsilon^{-\alpha}<\infty$.
Thus convergence of the truncated point processes implies convergence
of the corresponding finite point configurations. After adjoining the
origin, their convex hulls therefore converge in the Hausdorff metric.
Since the circumradius is continuous with respect to Hausdorff
convergence of compact sets, the continuous mapping theorem yields
\[
C\bigl(\operatorname{conv}(\{0\}\cup N_n^{(\varepsilon)})\bigr)
\xrightarrow{\mathcal D}
C\bigl(\operatorname{conv}(\{0\}\cup N^{(\varepsilon)})\bigr).
\]

It remains to control the points omitted by the truncation. To this end, put
\[
A_n^{(\varepsilon)}
:=
\{c_n^{-1}X_k:\|X_k\|>\varepsilon c_n,\ 1\le k\le n\}, \quad
B_n^{(\varepsilon)}
:=
\{c_n^{-1}X_k:\|X_k\|\le\varepsilon c_n,\ 1\le k\le n\}.
\]
Then $B_n^{(\varepsilon)}\subset B(0,\varepsilon)$. Let
$S_n^{(\varepsilon)} := A_n^{(\varepsilon)}\cup B_n^{(\varepsilon)}
= \{c_n^{-1}X_1,\ldots,c_n^{-1}X_n\}$
and
$T_n^{(\varepsilon)} := \{0\}\cup A_n^{(\varepsilon)}$.
Notice that $C(S_n^{(\varepsilon)})=c_n^{-1}R_n$
and $C(T_n^{(\varepsilon)}) = C\bigl(\operatorname{conv}(\{0\}\cup N_n^{(\varepsilon)})\bigr)$.

We first observe that $C(S_n^{(\varepsilon)}) \le C(T_n^{(\varepsilon)})+\varepsilon$.
Indeed, if $B(c,r)$ contains $T_n^{(\varepsilon)}$, then $0\in
B(c,r)$, and hence $\|c\|\le r$. For every
$b\in B_n^{(\varepsilon)}$, $\|b-c\| \le \|b\|+\|c\| \le \varepsilon+r$.
Thus $B(c,r+\varepsilon)$ contains $S_n^{(\varepsilon)}$, which proves
the assertion.

For the converse inequality, suppose that
$B_n^{(\varepsilon)}\neq\varnothing$ and choose
$b_0\in B_n^{(\varepsilon)}$. If $B(c,r)$ contains
$S_n^{(\varepsilon)}$, then
$\|c\| \le \|c-b_0\|+\|b_0\| \le r+\varepsilon$.
Hence $0\in B(c,r+\varepsilon)$, and therefore
$B(c,r+\varepsilon)$ contains $T_n^{(\varepsilon)}$. Consequently,
on the event $\{B_n^{(\varepsilon)}\neq\varnothing\}$,
$\big| C(S_n^{(\varepsilon)}) - C(T_n^{(\varepsilon)}) \big| \le \varepsilon$.

For every fixed $\varepsilon>0$,
$\mathbb P\bigl(B_n^{(\varepsilon)}=\varnothing\bigr)
= \mathbb P(\|X\|>\varepsilon c_n)^n \rightarrow 0$,
since $\varepsilon c_n\to\infty$ and hence
$\mathbb P(\|X\|>\varepsilon c_n)\to0$. It follows that, for every
$\delta>0$ and every $\varepsilon<\delta$,
\[
\limsup_{n\to\infty} \mathbb P\left( \left| c_n^{-1}R_n
- C\bigl(\operatorname{conv}(\{0\}\cup N_n^{(\varepsilon)})\bigr)
\right|>\delta \right) =0.
\]
Therefore, 
\[
\lim_{\varepsilon\downarrow 0} \limsup_{n\to\infty} \mathbb P\left(
\left| c_n^{-1}R_n -
C\bigl(\operatorname{conv}(\{0\}\cup N_n^{(\varepsilon)})\bigr) \right|>\delta
\right) =0.
\]

It remains to prove the corresponding assertion for the limiting
Poisson cloud. Put
\[
S:=\{0\}\cup\{P_k:k\ge1\},
\qquad
T^{(\varepsilon)}
:=
\{0\}\cup\{P_k:\|P_k\|>\varepsilon\}.
\]
Clearly $T^{(\varepsilon)}\subset S$, and therefore
$C(T^{(\varepsilon)})\le C(S)$.

Conversely, if $B(c,r)$ contains $T^{(\varepsilon)}$, then
$0\in B(c,r)$ and hence $\|c\|\le r$. Every point $P_k$ omitted by
the truncation satisfies $\|P_k\|\le\varepsilon$, and consequently
$\|P_k-c\| \le \|P_k\|+\|c\| \le\varepsilon+r$.
Thus $B(c,r+\varepsilon)$ contains $S$, and hence
$0 \le C(S)-C(T^{(\varepsilon)}) \le\varepsilon$ a.s.
Since the circumradius is unchanged by taking convex hulls, we have
$C(T^{(\varepsilon)})= 
C\bigl(\operatorname{conv}(\{0\}\cup N^{(\varepsilon)})\bigr)$.
It follows that, for every $\delta>0$,
\[
\lim_{\varepsilon\downarrow0}
\mathbb P\left(
\left|
C(\{0\}\cup\{P_k:k\ge1\})
-
C\bigl(\operatorname{conv}(\{0\}\cup N^{(\varepsilon)})\bigr)
\right|>\delta
\right)
=0.
\]
This completes the proof.
\end{proof}

We are now in a position to derive the Fr\'{e}chet limit theorem
for the enclosing radius. Lemma~\ref{lem:Frechet-truncation} shows that, after truncation away from the origin, the relevant geometric functionals are
continuous under vague convergence and that the omitted points
do not affect the enclosing radius asymptotically.

For the diameter $D_n$, limit theorems in the Fréchet regime were
obtained by Jammalamadaka and Janson \cite{jamsva15} and
by Henze and Lao \cite{henzelao}. The following result shows that
the enclosing radius $R_n$ admits a corresponding limit description in
terms of the same limiting Poisson cloud of extreme observations.
The resulting limit law is given by a functional of the limiting Poisson cloud
and typically does not admit a simple closed form.

Since the Poisson point process $N=\sum_{k\ge1}\varepsilon_{P_k}$
has infinitely many points accumulating at the origin almost surely,
the origin belongs almost surely to the closure of
$\{P_k:k\ge1\}$. Since the circumradius depends only on the
closure of a set, we therefore have
$C(\{P_k:k\ge1\}) = C(\{0\}\cup\{P_k:k\ge1\})$ a.s. 
Consequently, the origin need not appear explicitly in the limit
statement below.

\begin{theorem}\label{thm:Frechet-radius}
Assume that \eqref{eq-reg-var} holds. Let $P_1,P_2,\ldots$ denote the
points of a Poisson point process on $\mathbb{R}^d\setminus\{0\}$
with intensity measure $\nu$ given in \eqref{sec-4-intmeas}. Then
\[
        c_n^{-1}R_n
        \vertk
        \mathcal{C}\bigl(\{P_k:k\ge1\}\bigr).
\]
\end{theorem}

\smallskip
\noindent 
\begin{proof}
For $\varepsilon>0$, put $Y_n^{(\varepsilon)} := C\bigl(\operatorname{conv}(\{0\}\cup N_n^{(\varepsilon)})\bigr)$ and
$Y^{(\varepsilon)} := C\bigl(\operatorname{conv}(\{0\}\cup N^{(\varepsilon)})\bigr)$.
By Lemma~\ref{lem:Frechet-truncation}, we have 
$Y_n^{(\varepsilon)} \vertk Y^{(\varepsilon)}$
as $n\to\infty$ for each $\varepsilon >0$.
Moreover, for every $\delta>0$,
\[
\lim_{\varepsilon\downarrow0} \limsup_{n\to\infty} 
\mathbb P\left(\left| c_n^{-1}R_n-Y_n^{(\varepsilon)} \right|>\delta \right) =0,
\]
and
\[
\lim_{\varepsilon\downarrow0} \mathbb P\left(
\left| C(\{0\}\cup\{P_k:k\geq1\}) - Y^{(\varepsilon)} \right|>\delta \right) =0.
\]
Hence
$c_n^{-1}R_n \vertk C(\{0\}\cup\{P_k:k\geq1\})$ by the converging-together theorem.
Since the origin belongs almost surely to the closure of
$\{P_k:k\geq1\}$ and the circumradius is unchanged upon taking
closures, we have
$C(\{0\}\cup\{P_k:k\geq1\}) = C(\{P_k:k\geq1\})$ almost surely.
Therefore $c_n^{-1}R_n \vertk C(\{P_k:k\geq1\})$, 
which proves the assertion.
\end{proof}

\begin{corollary}
Suppose in addition that $1-F(s)\sim \gamma s^{-\alpha}$ as $s\to\infty$,
for some $\gamma >0$. Then $c_n\sim(\gamma n)^{1/\alpha}$ and therefore
\[
(\gamma n)^{-1/\alpha}R_n
\vertk
\mathcal{C}(\{P_k:k\ge1\}).
\]
\end{corollary}

In dimension $d = 2$ the smallest enclosing ball is the smallest
enclosing circle. It is well known that such a circle is determined
either by two points forming a diameter or by three points forming
the circumcircle of an acute triangle. Thus the support size of the
minimal enclosing circle is either $2$ or $3$ in the planar case.

\medskip
\noindent
To illustrate the geometry arising in the Fréchet regime,
we performed a small simulation study in dimension $d=2$.
We generated i.i.d.\ vectors in $\mathbb{R}^2$ with uniform angular
distribution and Pareto radial tail
$\mathbb{P}(\|X\| >r)=r^{-\alpha}$, $r\ge 1$,
for several values of the tail index $\alpha$.
For each simulated sample we computed the smallest enclosing circle
and recorded whether its support consisted of $k=2$ or $k=3$ points.

As a general notation, let $K_n$ denote the  number of sample
points supporting the smallest enclosing ball.
Under the standing assumption of absolute continuity, 
we have $K_n\in\{2,\ldots,d+1\}$ almost surely.
Thus $K_n\in\{2,3\}$ almost surely in the planar case.

For all simulation studies reported in this paper, we used the following
computational procedure. For each sample we put $m=\min\{50,n\}$ and begin
by applying Welzl's algorithm (see Welzl~\cite{wel91} and Fischer et al.~\cite{fi03}),
including the move-to-front heuristic, to the
$m$ points of largest norm. We then scan all $n$ points for violators, that
is, points lying outside the resulting ball, and add all violators to the
active set to which Welzl's algorithm is applied.

Whenever at least one violator is added, the radius of the recomputed smallest
enclosing ball increases strictly. Since the active set grows by at
least one point in every such round, the procedure terminates after
finitely many rounds. Once no violator remains, the ball obtained from
the active set is the smallest enclosing ball of the entire sample:
it contains every sample point and hence its radius is at least the
optimal radius for the full sample, while optimality for the active
subset implies the reverse inequality.

To verify the outcome of every replication, we use the fact that a ball
containing all sample points is the smallest enclosing ball if and only
if its centre lies in the convex hull of the sample points on its
boundary sphere; see Welzl~\cite{wel91} and Fischer et al.~\cite{fi03}.
We therefore check that (i) all $n$ sample points lie in the computed
ball, (ii) the $k$ support points returned by the algorithm lie on its
boundary sphere, are affinely independent, and the barycentric
coordinates of the centre with respect to these points are strictly
positive, and (iii) exactly $k$ sample points lie within relative
tolerance $10^{-10}$ of the boundary sphere. Conditions (i) and (ii)
imply that the computed ball is the smallest enclosing ball; for this
conclusion, non-negative barycentric coordinates would already suffice.
Strict positivity additionally implies that no proper subset of these
$k$ support points has the centre in its convex hull. Together with
condition (iii), this verifies that the support set is unique and has
size $k$.

Replications that fail one of these checks are recorded separately but
retained in the reported statistics. No replication failed any check
except in the Pareto model with $\alpha=1/2$, where the fraction of
flagged replications was below $5\times10^{-5}$. Thus, even if all
flagged replications were classified differently, the resulting change
in any reported support probability would be at most $5\times10^{-5}$,
which is negligible on the scale of the Monte Carlo error.

Table~\ref{tab:frechet-support} shows the estimated probabilities $\mathbb{P}(K_n=j)$,
$j=2,3$, for several values of $n$ and $\alpha$.

\begin{table}[h]
\centering
\begin{tabular}{c|r|cc}
$\alpha$ & $n$ & $\widehat{\mathbb{P}}(K_n=2)$ & $\widehat{\mathbb{P}}(K_n=3)$ \\
\hline
0.5 & 50   & 0.837 & 0.163 \\
0.5 & 200  & 0.837 & 0.163 \\
0.5 & 1000 & 0.838 & 0.163 \\
\hline
1   & 50   & 0.711 & 0.289 \\
1   & 200  & 0.712 & 0.288 \\
1   & 1000 & 0.711 & 0.289\\
\hline
2   & 50   & 0.547 & 0.454 \\
2   & 200  & 0.545 & 0.455 \\
2   & 1000 & 0.545 & 0.455 \\
\hline
4   & 50   & 0.375 & 0.623 \\
4   & 200  & 0.373 & 0.627 \\
4   & 1000 & 0.374 & 0.626 \\
\end{tabular}
\caption{Estimated probabilities of the support size of the minimal
enclosing circle in the planar Fréchet model ($10^6$ replications).}
\label{tab:frechet-support}
\end{table}

Two features are worth noting. First, for fixed $\alpha$, the estimated
probabilities are already remarkably stable across the values of $n$. This
behaviour is consistent with the scale invariance of the limiting Poisson
cloud that appears in the theorem.

Second, the balance between the cases $k=2$ and $k=3$ depends substantially
on the tail index $\alpha$. For heavier tails, corresponding to smaller
values of $\alpha$, the extreme radial observations tend to be more strongly
separated in magnitude, so that the enclosing circle is more often determined
by only two exceptionally distant points. As $\alpha$ increases, the largest
radii become less dominant relative to the remaining extreme observations,
and three-point configurations occur more frequently.

These observations admit a natural interpretation in view of
Theorem~\ref{thm:Frechet-radius}. After normalization, the extreme observations converge to a Poisson point process in $\mathbb R^2\setminus\{0\}$, and the
limiting enclosing circle is therefore the smallest circle containing
this Poisson cloud. Let $K$ denote the number of Poisson points that
determine this minimal enclosing circle. Then $K\in\{2,3\}$ almost
surely.

Here the supporting points are genuine points of the Poisson process:
since the process has only finitely many points outside every
neighbourhood of the origin, its only accumulation point is the origin,
whereas the minimal enclosing circle has positive radius almost surely.

The simulation results in Table~\ref{tab:frechet-support} are consistent with this
limiting Poisson geometry and suggest that both two- and three-point
support configurations persist asymptotically in the planar
Fr\'echet regime.

This behaviour contrasts with the corresponding bounded case
(for example when the sample is drawn from the unit disk), where the
asymptotic geometry of the enclosing circle is governed by points
accumulating near the boundary and the support structure exhibits a
different limiting pattern (see Section~\ref{secsuffcond}).

We also carried out a small simulation study in dimension $d=3$ for the same
Fr\'echet model. In this case the smallest enclosing ball is determined by
$k=2,3$, or $4$ sample points. As before, we generated i.i.d.\ random vectors
in $\mathbb{R}^3$ with uniform angular distribution on $\mathbb{S}^2$ and Pareto radial
tail $\mathbb{P}(\|X\|>r)=r^{-\alpha}$, $r\ge 1$,
and for each sample we computed the smallest enclosing ball and recorded the
support size $K_n$.

Table~\ref{tabsec3-2} shows the estimated probabilities $\widehat{\mathbb{P}}(K_n=j)$,
$j=2,3,4$, for several values of $n$ and $\alpha$.

\begin{table}[h]
\centering
\begin{tabular}{ccccc}
$\alpha$ & $n$ & $\widehat{\mathbb{P}}(K_n=2)$ & $\widehat{\mathbb{P}}(K_n=3)$ & $\widehat{\mathbb{P}}(K_n=4)$\\
\hline
0.5 & 50   & 0.768 & 0.217 & 0.016 \\
0.5 & 200  & 0.768 & 0.217 & 0.016 \\
0.5 & 1000 & 0.767 & 0.217 & 0.016 \\
1   & 50   & 0.596 & 0.353 & 0.051 \\
1   & 200  & 0.596 & 0.353 & 0.051 \\
1   & 1000 & 0.587 & 0.352 & 0.051 \\
2   & 50   & 0.390 & 0.472 & 0.138 \\
2   & 200  & 0.389 & 0.473 & 0.138 \\
2   & 1000 & 0.389 & 0.473 & 0.138 \\
4   & 50   & 0.206 & 0.502 & 0.292 \\
4   & 200  & 0.205 & 0.503 & 0.292 \\
4   & 1000 & 0.205 & 0.503 & 0.292 \\
\hline
\end{tabular}

\smallskip
\caption{Estimated probabilities of the support size of the minimal
enclosing ball in the three-dimensional Fr\'echet model ($10^6$ replications).}
\label{tabsec3-2}
\end{table}

As in the planar case, the estimated probabilities are remarkably stable
across the three sample sizes for each fixed value of $\alpha$. This
stability is consistent with the limiting Poisson-cloud description in
Theorem~\ref{thm:Frechet-radius} and suggests that the support-size distribution approaches its
limiting regime already for rather moderate sample sizes. The
three-dimensional picture is nevertheless more varied than in the planar
case. For heavy tails (small $\alpha$), the two-point configuration
dominates, reflecting the stronger radial separation among the most extreme
observations. As $\alpha$ increases, several extreme points tend to have
more comparable magnitudes, and the three-point and four-point
configurations become markedly more frequent. For $\alpha=4$, all three
support sizes occur with substantial probability. Thus, in contrast to the
planar bounded-support case (see Section~\ref{secsuffcond}), the Fr\'echet regime in
dimension three does not appear to favour a single support size
overwhelmingly; rather, the geometry of the minimal enclosing ball reflects
a genuine competition between two-, three-, and four-point extremal
configurations.

\section{The Weibull case}\label{secweib}
In this section we consider i.i.d.\ spherically symmetric random vectors
$X_1,X_2,\ldots$ in $\mathbb{R}^d$, $d\geq 2$, whose distribution is
supported on the closed unit ball $\mathbb{B}^d=\{x\in\mathbb{R}^d:\|x\|\leq 1\}$.
As before, we write $X_i=\rho_i U_i$, $i\geq 1$, where
the random directions $U_i$ are uniformly distributed on
$\mathbb{S}^{d-1}$ and independent of the radial variables $\rho_i = \|X_i\|$.
We assume that, for some constants $C>0$ and $\beta>-1$,
\begin{equation} \label{eq:boundary-tail}
  \mathbb{P}(\rho_1>1-s)
  \sim C s^{\beta+1},
  \qquad s\downarrow 0.
  \end{equation}

Taken together, Sections~\ref{secgumbelkotz}, \ref{secfrechet} and~\ref{secweib} address the enclosing-radius problem
in the three classical extreme-value regimes: Gumbel, Fr\'echet, and
Weibull. In contrast to the Gumbel and Fr\'echet regimes, the relevant
extreme observations in the bounded-support case accumulate near the
boundary of the support rather than appearing as isolated radial
outliers. The natural extremal object is therefore a point process
describing simultaneously the directions of sample points and their
rescaled distances from the boundary.

\subsection{Boundary point process}
Put $a_n:=n^{1/(\beta+1)}$ 
and define the boundary depths $Y_{ni}:=a_n(1-\rho_i)$ for $1\leq i\leq n$.
We consider the point process
\begin{equation}\label{eq:defpnn}
N_n:=\sum_{i=1}^n \varepsilon_{(U_i,Y_{ni})}
\end{equation}
on $\mathbb{S}^{d-1}\times[0,\infty)$. The following result identifies the limiting boundary point process.

\medskip

\begin{theorem}\label{thm:boundary-pp}
Let $X_1,X_2,\ldots$ be i.i.d.\ spherically symmetric random vectors
in $\mathbb{B}^d$, $d\geq 2$, and suppose that \eqref{eq:boundary-tail} holds.
Then
\[
N_n \xrightarrow{\mathcal D} N
\]
in $M_P(\mathbb{S}^{d-1}\times[0,\infty))$, equipped with the vague
topology, where $N$ is a Poisson point process with intensity measure
\begin{equation}
  \mu(\D u,\D t)
  =
  C(\beta+1)t^\beta\,\D t\,\sigma_{d-1}(\D u).
  \label{eq:boundary-intensity}
\end{equation}
Here $\sigma_{d-1}$ denotes the uniform probability measure on
$\mathbb{S}^{d-1}$.
\end{theorem}

\smallskip
\begin{proof}
We verify convergence of Laplace functionals. To this end, let 
$f:\mathbb{S}^{d-1}\times[0,\infty)\rightarrow[0,\infty)
$
be continuous with compact support. The Laplace functional of the limiting Poisson process $N$ is
\[
\mathbb{E}\bigl[\E^{-N(f)}\bigr]
\ = \ 
\exp\left\{   -   \int_{\mathbb{S}^{d-1}}
  \int_0^\infty   \bigl(1-\E^{-f(u,t)}\bigr)
  C(\beta+1)t^\beta\,\D t\,   \sigma_{d-1}(\D u) \right\}.
\]
Since the pairs
$(U_i,\rho_i)$ are i.i.d., we have
\[
\mathbb{E}\bigl[\E^{-N_n(f)}\bigr]
= \left(   \mathbb{E}   \left[    \E^{-f(U_1,Y_{n1})}   \right]\right)^n
=\left(  1-   \mathbb{E}   \left[    1-\E^{-f(U_1,Y_{n1})}   \right]
\right)^n.
\]
Thus, in view of $(1-t_n)^n\to \E^{-t}$ whenever
$nt_n\to t$, it suffices to show that
\begin{equation}
 n\,\mathbb{E}  \left[    1-\E ^{-f(U_1,Y_{n1})}
 \right]  \rightarrow  \int_{\mathbb{S}^{d-1}}
 \int_0^\infty  \bigl(1-\E^{-f(u,t)}\bigr)
 C(\beta+1)t^\beta\,\D t\,  \sigma_{d-1}(\D u).
 \label{eq:laplace-target}
\end{equation}

Define finite measures $\nu_n$ on $[0,\infty)$ by
$\nu_n(B):= n\,\mathbb{P}(Y_{n1}\in B)$.
Using $a_n^{\beta+1}=n$ and
\eqref{eq:boundary-tail}, we obtain
\[
\nu_n([0,t])
= n\,\mathbb{P}
\left(   \rho_1\geq 1-\frac{t}{a_n}
\right) \rightarrow C t^{\beta+1}
\]
for each $t >0$. Consequently,
$\nu_n \xrightarrow{v} \nu$
on $[0,\infty)$, where
$\nu(\D t)=C(\beta+1)t^\beta\,\D t$.

By independence of $U_1$ and $\rho_1$,
\[
n\,\mathbb{E}
 \left[    1-\E ^{-f(U_1,Y_{n1})}  \right]
 = \int_{\mathbb{S}^{d-1}}
\left(   \int_0^\infty   \bigl(1-\E ^{-f(u,t)}\bigr)\,\nu_n(\D t)
\right) \sigma_{d-1}(\D u).
\]
For fixed $u\in\mathbb{S}^{d-1}$, the function
$g_u(t):=1-\E^{-f(u,t)} $
is continuous and has compact support. Hence, by vague convergence,
\[
\int_0^\infty g_u(t)\,\nu_n(\D t)
\rightarrow
\int_0^\infty
g_u(t)C(\beta+1)t^\beta\,\D t.
\]

It remains only to justify integration with respect to $u$. Since $f$
has compact support, there exists $T<\infty$ such that
$f(u,t)=0$ for all $u\in\mathbb{S}^{d-1}$ and $t>T$. Thus
\[
0\leq \int_0^\infty g_u(t)\,\nu_n(\D t)
\leq \nu_n([0,T]),
\]
uniformly in $u$. Moreover,
$\nu_n([0,T])\rightarrow CT^{\beta+1},
$
and hence the right-hand side is bounded uniformly in $n$ for all
sufficiently large $n$. Dominated convergence therefore yields
\eqref{eq:laplace-target}, as was to be shown. 
\end{proof}

\medskip

\begin{remark}\label{rem:uniform-ball}
For the uniform distribution on $\mathbb{B}^d$,
$\mathbb{P}(\rho_1>1-s) = 1-(1-s)^d \sim ds$ as $s\downarrow 0$ .
Thus \eqref{eq:boundary-tail} holds with $\beta=0$ and $C=d$, and
$a_n=n$. In this case Theorem~\ref{thm:boundary-pp} yields a Poisson
point process on $\mathbb{S}^{d-1}\times[0,\infty)$ with intensity
$\mu(\D u,\D t)=d\,\sigma_{d-1}(\D u)\,\D t$.
\end{remark}

Theorem~\ref{thm:boundary-pp} identifies the sample points that are relevant on the
natural boundary scale. The next step is to relate this boundary point
process to the smallest enclosing ball. To this end, we first derive a
uniform first-order approximation of the condition that a ball contains
a point lying close to the boundary of $B^d$.

\subsection{Geometric linearization and a variational representation}\label{subseclinear}

We next relate the boundary point process to the smallest enclosing
ball. For $u\in\mathbb S^{d-1}$, $t\geq 0$, and $z\in\mathbb R^d$, put
\[ 
h_n(u,t;z) := a_n\left( 1- \left\| \left(1-\frac{t}{a_n}\right)u-\frac{z}{a_n}
\right\| \right).
\]
The following lemma gives the first-order approximation of this quantity on the natural boundary scale.

\medskip

\begin{lemma}\label{lem:weibull-linearization}
For every $M,L>0$,
\[
\sup_{\substack{u\in\mathbb S^{d-1}\\
                 0\leq t\leq M\\
                 \|z\|\leq L}}
\left| h_n(u,t;z)-\bigl(t+ z\cdot u\bigr)
\right| \rightarrow 0
\]
as $n\to\infty$. Moreover, for all $u\in\mathbb S^{d-1}$, $0 \le t \le a_n$, and
$z\in\mathbb R^d$, we have 
\begin{equation}\label{eq:weibull-bounds}
 t-\|z\|  \leq h_n(u,t;z)  \leq t+ z\cdot u .
 \end{equation}
\end{lemma}

\smallskip
\begin{proof}
For fixed $M,L>0$, write
\[
v_n= \left(1-\frac{t}{a_n}\right)u-\frac{z}{a_n}.
\]
Then
\[
\|v_n\|^2 = 1-\frac{2(t+ z\cdot u)}{a_n}
+\frac{\|tu+z\|^2}{a_n^2}.
\]
Uniformly for $u\in\mathbb S^{d-1}$, $0\leq t\leq M$, and
$\|z\|\leq L$, Taylor expansion of the square root yields
\[
\|v_n\| = 1-\frac{t+ z\cdot u}{a_n} +O(a_n^{-2}).
\]
Consequently, $h_n(u,t;z) = t+ z\cdot u+O(a_n^{-1})$,
uniformly over the stated ranges, proving the first assertion.
For the lower bound in \eqref{eq:weibull-bounds}, the triangle
inequality gives
\[
\left\| \left(1-\frac{t}{a_n}\right)u-\frac{z}{a_n} \right\|
\leq 1-\frac{t}{a_n}+\frac{\|z\|}{a_n},
\]
and hence $h_n(u,t;z)\geq t-\|z\|$.
For the upper bound, using $\|x\|\geq x\cdot u$ for
$\|u\|=1$, and putting 
\[
x = \left((1- \frac{t}{a_n}\right) u- \frac{z}{a_n},
\]
we obtain
\[
\left\| \left(1-\frac{t}{a_n}\right)u-\frac{z}{a_n} \right\|
\geq 1-\frac{t}{a_n}-\frac{z\cdot  u}{a_n},
\]
which yields
$h_n(u,t;z)\leq t+ z\cdot u$.
\end{proof}

For the sample points, recall that
$X_i= (1-Y_{ni}/a_n)U_i$.
Define
\begin{equation}\label{eq:defhnz}
H_n(z) := \min_{1\leq i\leq n} h_n(U_i,Y_{ni};z), \qquad z\in\mathbb R^d,
\end{equation}
and its linearized counterpart
\begin{equation}\label{eq:defgnz}
G_n(z) := \min_{1\leq i\leq n} \bigl(Y_{ni}+ z\cdot U_i\bigr).
\end{equation}

The definition of the smallest enclosing radius gives the exact
variational identity
\begin{equation}\label{eq:exact-weibull-variational}
a_n(1-R_n)
=
\sup_{z\in\mathbb R^d} H_n(z).
\end{equation}
Indeed,
\[
a_n(1-R_n) =
a_n\left( 1-\inf_{c\in\mathbb R^d}        \max_{1\leq i\leq n}\|X_i-c\| \right)
= \sup_{z\in\mathbb R^d} \min_{1\leq i\leq n} a_n\left(
1-\left\|X_i-\frac{z}{a_n}\right\| \right).
\]
In particular, if $C_n$ denotes the centre of the smallest enclosing
ball, then $a_nC_n$ is a maximizer of $H_n$.


\subsection{Tightness of the centres and truncation}\label{subsec-centres}

To pass from the exact variational representation \eqref{eq:exact-weibull-variational}
to a limit theorem, two compactness issues have to be addressed. We first
control the optimizing centres on the scale $a_n^{-1}$ and then show that
points sufficiently far from the boundary are asymptotically irrelevant.

We first show that the centres of the smallest enclosing balls are tight on the natural boundary scale. Recall the functions $H_n$ and $G_n$ introduced in Subsection~\ref{subseclinear}.

\medskip

\begin{lemma}\label{lem:weibull-centre-tightness}
Let $C_n$ denote the centre of the smallest enclosing ball of
$X_1,\ldots,X_n$. Under \eqref{eq:boundary-tail},
\[
a_n C_n=O_{\mathbb P}(1).
\]
Moreover, for every $\varepsilon>0$ there exists $L>0$ such that,
for all sufficiently large $n$,
\[
\mathbb P\left(
\sup_{z\in\mathbb R^d}H_n(z)
=
\sup_{\|z\|\leq L}H_n(z)
\right)\geq 1-\varepsilon.
\]
The same assertion holds with $H_n$ replaced by $G_n$.
\end{lemma}

\smallskip
\begin{proof}
Fix $\varepsilon>0$. Choose $\delta\in(0,1/\sqrt{2})$ and points
$v_1,\ldots,v_m\in\mathbb S^{d-1}$ such that the spherical neighbourhoods
$B(v_j,\delta) = \{u\in\mathbb S^{d-1}:\|u-v_j\|<\delta\}$ for $j=1,\ldots,m$
cover $\mathbb S^{d-1}$.  For $M>0$ and $j=1,\ldots,m$, define
$A_{n,j} = \left\{ N_n\bigl(B(-v_j,\delta)\times[0,M]\bigr)\geq 1 \right\}$.
Thus, on $A_{n,j}$ there exists at least one sample point whose
direction lies in $B(-v_j,\delta)$ and whose boundary depth is at most
$M$.  Writing $q_j:=\sigma_{d-1}\bigl(B(-v_j,\delta)\bigr)>0$,
we have
\[
\mu\bigl(B(-v_j,\delta)\times[0,M]\bigr) = Cq_jM^{\beta+1}.
\]
By Theorem~\ref{thm:boundary-pp}, it follows that 
$N_n\bigl(B(-v_j,\delta)\times[0,M]\bigr)
\vertk N\bigl(B(-v_j,\delta)\times[0,M]\bigr)$,
since the boundary of this set carries no mass under $N$ almost surely.
The limiting random variable is Poisson with mean $Cq_jM^{\beta+1}$. Consequently,
$\mathbb P(A_{n,j}^c) \rightarrow \exp\{-Cq_jM^{\beta+1}\}$.

Put $E_n:=\bigcap_{j=1}^m A_{n,j}$. Since there are only finitely many caps, $M$ can be chosen so large that
\begin{equation}\label{eq:liminfen}
\liminf_{n\to\infty}\mathbb P(E_n)\geq1-\frac{\varepsilon}{2}
\end{equation}
and thus $\mathbb P(E_n)\ge 1-\varepsilon$ for sufficiently large $n$.

Now, put $\eta:=1-2\delta^2$. Note that $\eta >0$ since $\delta < 1/\sqrt{2}$.
On $E_n$, the following property holds. For every
$v\in\mathbb S^{d-1}$ choose $j$ such that
$v\in B(v_j,\delta)$. Then there exists an index $i=i(v)$ such that
$Y_{ni}\leq M$ and $U_i\in B(-v_j,\delta)$. Consequently,
$\|U_i+v\| \leq \|U_i+v_j\|+\|v_j-v\| <2\delta$. Since $\|U_i\|=\|v\|=1$, we have
$2+2U_i\cdot v=\|U_i+v\|^2<4\delta^2$ and therefore
$U_i\cdot v<2\delta^2-1=-\eta$.

Now let $z\neq0$ and put $v=z/\|z\|$. By the upper bound in
\eqref{eq:weibull-bounds},
\[
H_n(z) \leq h_n(U_i,Y_{ni};z) \leq Y_{ni}+z\cdot U_i < M-\eta\|z\|.
\]
On the other hand,
$H_n(0)=\min_{1\leq i\leq n}Y_{ni}\ge 0$.
Thus, on $E_n$, $H_n(z)<H_n(0)$ whenever $\|z\|>M/\eta$.
It follows that every maximizer of $H_n$ lies in the ball
$\{z:\|z\|\leq M/\eta\}$ on $E_n$. Since $a_nC_n$ is a maximizer
of $H_n$, this and \eqref{eq:liminfen} prove
$a_nC_n=O_{\mathbb P}(1)$ and the compactness assertion for $H_n$.

The same argument applies to $G_n$, since on $E_n$
\[
G_n(z)
\leq Y_{ni}+z\cdot U_i
\leq M-\eta\|z\|,
\qquad
G_n(0)=\min_{1\leq i\leq n}Y_{ni}\geq0.
\]
\end{proof}

We next show that, once the optimization is restricted to a compact
set, points lying sufficiently far from the boundary cannot affect the
variational problem. To this end, let $M >0$ and define
\begin{eqnarray}\nonumber
H_{n,M}(z) & := & \min_{\substack{1\leq i\leq n\\ Y_{ni}\leq M}} h_n(U_i,Y_{ni};z),\\ \label{eq:defgnmz} 
G_{n,M}(z) & := &
\min_{\substack{1\leq i\leq n\\ Y_{ni}\leq M}}
\bigl(Y_{ni}+z\cdot U_i\bigr),
\end{eqnarray}
where the minimum over an empty set is interpreted as $+\infty$.

\medskip

\begin{lemma}\label{lem:weibull-truncation}
For every $L>0$,
\[
\lim_{M\to\infty}\liminf_{n\to\infty}
\mathbb P\Big(H_n(z)=H_{n,M}(z)\ \text{for all }\|z\|\leq L \Big)=1,
\]
and
\[
\lim_{M\to\infty}\liminf_{n\to\infty}
\mathbb P\Big(G_n(z)=G_{n,M}(z)\ \text{for all }\|z\|\leq L\Big)=1.
\]
\end{lemma}

\smallskip
\begin{proof}
Fix $L>0$ and let $M>4L$. Consider the event
\begin{equation}\label{eq:defbnm}
B_{n,M} := \left\{ \min_{1\leq i\leq n}Y_{ni}\leq \frac{M}{2} \right\}.
\end{equation}
On $B_{n,M}$, choose an index $i$ such that $Y_{ni}\leq M/2$.
For $\|z\|\leq L$, the upper bound in \eqref{eq:weibull-bounds}
gives
$h_n(U_i,Y_{ni};z) \leq Y_{ni}+z\cdot U_i \le M/2 +L$. 
On the other hand, if $Y_{nj}>M$, then the lower bound in
\eqref{eq:weibull-bounds} yields
$h_n(U_j,Y_{nj};z) \ge  Y_{nj}-\|z\| > M-L$.
Since $M>4L$, we have $M/2+L<M-L$. Hence, for every $\|z\|\leq L$ and every index $j$ with
$Y_{nj}>M$,
\[
h_n(U_i,Y_{ni};z) \leq \frac{M}{2}+L  <M-L < h_n(U_j,Y_{nj};z).
\]
Thus no index $j$ with $Y_{nj}>M$ can attain the minimum
defining $H_n(z)$, uniformly over $\|z\|\leq L$.
 Hence, on $B_{n,M}$, we have
$H_n(z)=H_{n,M}(z)$ if $\|z\|\le  L$. 
Using $Y_{nj}+z\cdot U_j\geq Y_{nj}-L$ the same argument shows that
$G_n(z)=G_{n,M}(z)$ if $\|z\|\leq L$,
on the same event. By Theorem~\ref{thm:boundary-pp},
\[
N_n\left(
\mathbb S^{d-1}\times \left[0,\frac{M}{2}\right] \right)
\vertk N\left( \mathbb S^{d-1}\times \left[0,\frac{M}{2}\right] \right),
\]
and the limiting random variable is Poisson with mean $C(M/2)^{\beta+1}$.  Therefore,
\begin{equation}\label{eq:convpbn}
\mathbb P(B_{n,M}^c) \rightarrow 
\exp\left\{-C\left(\frac{M}{2}\right)^{\beta+1}\right\}.
\end{equation}
Letting $M\to\infty$ proves both assertions.
\end{proof}


\subsection{Continuity of the truncated variational functional}\label{subseccontrunc}

Lemmas~\ref{lem:weibull-centre-tightness} and~\ref{lem:weibull-truncation} show that, with arbitrarily large probability, both the optimization over the centre variable and the minimization over the boundary depths may be restricted to compact sets. For fixed truncation parameters, the variational problem therefore depends only on finitely many points of the limiting Poisson process. This allows us to establish the continuity needed for an application of the continuous mapping theorem.

For $M,L>0$ and a locally finite point measure
\[
m=\sum_j \varepsilon_{(u_j,t_j)}
\]
on $\mathbb S^{d-1}\times[0,\infty)$, define
\begin{equation}\label{eq:defphinz}
\Phi_{M,L}(m) := \sup_{\|z\|\leq L} \inf_{\substack{(u,t)\in m\\ t\leq M}} (t+z\cdot u) \, \in \overline{\mathbb{R}} := [-\infty,\infty],
\end{equation}
where $\overline{\mathbb{R}}$ is equipped with its usual topology, and 
the infimum over an empty set is interpreted as $+\infty$.

\medskip

\begin{lemma}\label{lemma-7-6}
For every $M,L>0$, the mapping $m\mapsto \Phi_{M,L}(m)$ is almost surely continuous at the Poisson point process $N$ from Theorem~\ref{thm:boundary-pp} with respect to vague convergence.
\end{lemma}

\smallskip

\begin{proof}
Fix $M,L>0$. Since $\mu\bigl(\mathbb S^{d-1}\times[0,M]\bigr) = CM^{\beta+1}<\infty$,
the restriction of $N$ to $\mathbb S^{d-1}\times[0,M]$ contains only finitely many points
almost surely. Moreover,
$N\bigl(\mathbb S^{d-1}\times\{M\}\bigr)=0$ almost surely, because the intensity measure $\mu$ has no atom at $M$.
Putting $K_M := N\bigl(\mathbb S^{d-1}\times[0,M]\bigr)$, 
we work on the event of probability one on which $K_M<\infty$
and $N(\mathbb S^{d-1}\times\{M\})=0$.
Let $(m_n)_{n\geq1}$ be a sequence of locally finite point
measures such that $m_n\to N$ vaguely. Since
$\mathbb S^{d-1}\times[0,M]$ is a compact continuity set of $N$, we have
$m_n\bigl(\mathbb S^{d-1}\times[0,M]\bigr)=K_M$
for all sufficiently large $n$.

If $K_M=0$, then the corresponding restrictions of $m_n$ are
empty for all sufficiently large $n$, and hence
$\Phi_{M,L}(m_n)=\Phi_{M,L}(N)=+\infty$ eventually.

Suppose now that $K_M\geq1$. The points of $N$ in
$\mathbb S^{d-1}\times[0,M]$ may be written as
\[
(u_1,t_1),\ldots,(u_{K_M},t_{K_M}).
\]
By Proposition 3.13 of Resnick \cite{resn}, for all sufficiently large $n$, the points of $m_n$ in the same set can be labelled as
\[
(u_1^{(n)},t_1^{(n)}),\ldots, (u_{K_M}^{(n)},t_{K_M}^{(n)})
\]
in such a way that $(u_j^{(n)},t_j^{(n)}) \rightarrow (u_j,t_j)$ for $j=1,\ldots,K_M$.
For $\|z\|\leq L$, define
\[
F_n(z) := \min_{1\leq j\leq K_M} \bigl(t_j^{(n)}+z\cdot u_j^{(n)}\bigr), \qquad
F(z)  := \min_{1\leq j\leq K_M} (t_j+z\cdot u_j).
\]
Then, for every $\|z\|\leq L$,
\[
|F_n(z)-F(z)| \leq \max_{1\leq j\leq K_M} \left( |t_j^{(n)}-t_j|
+ L\|u_j^{(n)}-u_j\| \right).
\]
The right-hand side does not depend on $z$ and converges to zero.
Consequently,
\[
\sup_{\|z\|\leq L}|F_ n(z)-F(z)| \rightarrow 0.
\]
It follows that
\[
\bigl|\Phi_{M,L}(m_n)-\Phi_{M,L}(N)\bigr|
= \left| \sup_{\|z\|\leq L}F_n(z) - \sup_{\|z\|\leq L}F(z) \right|  
\leq \sup_{\|z\|\leq L}|F_n(z)-F(z)| \rightarrow 0.
\]
This proves the assertion.
\end{proof}

Consequently, Theorem~\ref{thm:boundary-pp} and the continuous mapping theorem yield,
for every fixed $M,L>0$,
\[
\sup_{\|z\|\leq L}G_{n,M}(z)
= \Phi_{M,L}(N_n) \ \vertk \Phi_{M,L}(N),
\]
where $N_n$ is defined in \eqref{eq:defpnn}.


\subsection{The limiting Poisson functional}\label{subseclpf}

The preceding results yield convergence of the variational problem when both the boundary depths and the centre variable are restricted to compact sets. We now remove these truncations and introduce the limiting functional appearing in the limit theorem.
To this end, put
\begin{equation}\label{eq:defghochm}
G(z) := \inf_{(u,t)\in N}(t+z\cdot u), \qquad
G^{(M)}(z) := \inf_{\substack{(u,t)\in N\\ t\leq M}}
(t+z\cdot u), \quad z \in \mathbb{R}^d, 
\end{equation}
where $M >0$ and, as before, the infimum over an empty set is interpreted as
$+\infty$. For $L>0$, define
\begin{equation}\label{eq:defphiln}
\Phi_L(N) := \sup_{\|z\|\leq L}G(z), 
\end{equation}
and finally
\[
\Phi(N) := \sup_{z\in\mathbb R^d}G(z)
= \sup_{z\in\mathbb R^d} \inf_{(u,t)\in N}(t+z\cdot u).
\]

\medskip

\begin{lemma}\label{lemma-7-7}
For every $L>0$,
\[
\lim_{M\to\infty} \mathbb P\Big( G(z)=G^{(M)}(z)\ \text{for all }\|z\|\leq L \Big) =1.
\]
Moreover, $G(z)\rightarrow-\infty$ as $\|z\|\to\infty$ almost surely. In particular, $\Phi(N)$ is finite almost surely, and the supremum defining $\Phi(N)$ is attained.
\end{lemma}

\smallskip
\begin{proof} 
Fix $L>0$ and let $M>4L$. Consider the event
\[
B_M :=  \left\{N\left( \mathbb S^{d-1}\times \left[0,\frac M2\right] \right)\geq 1 \right\}.
\]
On $B_M$, choose $(u,t)\in N\cap \bigl(\mathbb S^{d-1}\times[0,M/2]\bigr)$.
Then $t+z\cdot u\leq M/2 + L$  for every $\|z\|\leq L$.
On the other hand, for any $(u',t')\in N$ with $t'>M$, we have
$t'+z\cdot u'\geq t'-L>M-L$. Since $M>4L$, it follows that 
$M/2 +L<M-L$.
Hence no point with boundary depth larger than $M$ can attain the
infimum defining $G(z)$, simultaneously for all $\|z\|\leq L$.
Thus, on $B_M$, we have $G(z)=G^{(M)}(z)$ if $\|z\|\leq L$. Moreover,
\[
\mathbb P(B_M^c) = \exp\left\{ -C\left(\frac M2\right)^{\beta+1} \right\},
\]
which tends to zero as $M\to\infty$.

To prove the second assertion, choose $\delta\in(0,1/\sqrt2)$ and points
$v_1,\ldots,v_m\in\mathbb S^{d-1}$ such that $B(v_1,\delta),\ldots,B(v_m,\delta)$ cover
$\mathbb S^{d-1}$. Furthermore, put $\eta:=1-2\delta^2>0$.
For every $j$, the set $B(-v_j,\delta)\times[0,\infty)$ has infinite $\mu$-measure. Hence, almost surely, it contains at
least one point $(u_j,t_j)$ of $N$. Since there are only finitely many caps, it follows that 
\[
T:=\max_{1\leq j\leq m}t_j<\infty
\qquad\text{almost surely}.
\]
Now, let $z\neq0$ and write $v=z/\|z\|$. Choose $j$ such that
$v\in B(v_j,\delta)$. As in the proof of Lemma~\ref{lem:weibull-centre-tightness},
$u_j\cdot v<-\eta$. Consequently,
$G(z) \leq t_j+z\cdot u_j \leq T-\eta\|z\|$.
Thus $G(z)\rightarrow-\infty$ as $\|z\|\to\infty$ almost surely.

Since $t\geq0$ for every point $(u,t)\in N$, we have
$G(z)\geq-\|z\|>-\infty$ for every $z\in\mathbb R^d$. On the other
hand, since $N$ is non-empty almost surely, $G(z)<\infty$ for every
$z$. Thus $G$ is real-valued on $\mathbb R^d$ almost surely.
Moreover, as an infimum of affine functions, $G$ is upper
semicontinuous. Together with the preceding coercivity this shows
that its supremum is finite and attained.
\end{proof}

We now combine the continuity result from Subsection~\ref{subseccontrunc} with the
truncation arguments for the empirical and limiting functionals. This first yields convergence of the linearized variational problem on every fixed compact set in the centre variable. The uniform linearization from Lemma~\ref{lem:weibull-linearization} then shows that the exact and linearized functionals are asymptotically equivalent on such compact sets.

Recall $G_n(z)$ from \eqref{eq:defgnz}, $H_n(z)$ from \eqref{eq:defhnz}, and $\Phi_L(N)$ from \eqref{eq:defphiln}. 

\medskip

\begin{lemma}\label{lemma-7-8}
For every $L>0$,
\[
\sup_{\|z\|\leq L}G_n(z) \vertk \Phi_L(N).
\]
Moreover,
\[
\sup_{\|z\|\leq L} |H_n(z)-G_n(z)| \stk 0.
\]
Consequently,
\[
\sup_{\|z\|\leq L}H_n(z) \vertk \Phi_L(N).
\]
\end{lemma}

\smallskip

\begin{proof}
Fix $L>0$. We first prove the convergence of the linearized
functional. By Theorem~\ref{thm:boundary-pp}, Lemma~\ref{lemma-7-6}, and the continuous mapping theorem,
\[
\sup_{\|z\|\leq L}G_{n,M}(z) = \Phi_{M,L}(N_n) \vertk \Phi_{M,L}(N)
\]
for every fixed $M>0$. By Lemma~\ref{lem:weibull-truncation},
\[
\lim_{M\to\infty}\limsup_{n\to\infty}
\mathbb P\left(\sup_{\|z\|\leq L}G_n(z) \neq \sup_{\|z\|\leq L}G_{n,M}(z) \right) =0.
\]
On the other hand, by \eqref{eq:defphinz} and \eqref{eq:defghochm}, 
\[
\Phi_{M,L}(N) = \sup_{\|z\|\leq L}G^{(M)}(z), \qquad
\Phi_L(N) = \sup_{\|z\|\leq L}G(z).
\]
Hence $\bigl\{ G(z)=G^{(M)}(z) \text{ for all }\|z\|\leq L \bigr\}
\subseteq \bigl\{ \Phi_{M,L}(N)=\Phi_L(N) \bigr\}$.
It therefore follows from Lemma~\ref{lemma-7-7} that
\[
\lim_{M\to\infty} \mathbb P\left(\Phi_{M,L}(N)\neq\Phi_L(N)\right) =0.
\]
We are therefore in the setting of the converging-together theorem: for every fixed $M>0$, 
\[
\sup_{\|z\|\leq L}G_{n,M}(z) \vertk \Phi_{M,L}(N),
\]
whereas
\[
\lim_{M\to\infty}\limsup_{n\to\infty} \mathbb P\left(
\sup_{\|z\|\leq L}G_n(z) \neq \sup_{\|z\|\leq L}G_{n,M}(z) \right)=0 
\]
and
\[
\lim_{M\to\infty} \mathbb P\left( \Phi_{M,L}(N)\neq\Phi_L(N) \right)=0.
\]
Hence the converging-together theorem implies
\begin{equation}\label{eq:firstass}
\sup_{\|z\|\leq L}G_n(z) \vertk \Phi_L(N).
\end{equation}

We next compare the exact and the linearized functionals. For
$M,L>0$, define
\[
\Delta_n(M,L) := \sup_{\substack{ u\in\mathbb S^{d-1}\\ 
0\leq t\leq M\\ \|z\|\leq L}} \left| h_n(u,t;z)-(t+z\cdot u) \right|.
\]
By Lemma~\ref{lem:weibull-linearization},
$\Delta_n(M,L)\rightarrow 0$ for every fixed $M,L>0$. Let $M>4L$ and recall the event
$B_{n,M}$ from \eqref{eq:defbnm}. 
On $B_{n,M}$, both truncations are exact, that is, we have
$H_n(z)=H_{n,M}(z)$ and $G_n(z)=G_{n,M}(z)$ for $\|z\|\leq L$.
Moreover, since $B_{n,M}$ implies that $Y_{ni}\leq M/2<M$
for at least one index $i$, the truncated index set
$\{i:Y_{ni}\leq M\}$ is nonempty.  Therefore, for
$\|z\|\leq L$,
\[
\begin{split}
|H_n(z)-G_n(z)|
&=
\left|
\min_{\substack{1\leq i\leq n\\ Y_{ni}\leq M}}
h_n(U_i,Y_{ni};z)
-
\min_{\substack{1\leq i\leq n\\ Y_{ni}\leq M}}
(Y_{ni}+z\cdot U_i)
\right|\\
&\leq
\Delta_n(M,L).
\end{split}
\]
Consequently, for every $\varepsilon>0$, since
$\Delta_n(M,L)\to0$, we have for all sufficiently large $n$,
\[
\mathbb P\left(
\sup_{\|z\|\leq L}|H_n(z)-G_n(z)|>\varepsilon
\right)
\leq
\mathbb P(B_{n,M}^c).
\]
In view of \eqref{eq:convpbn}, it follows that
\[
\limsup_{n\to\infty} \mathbb P\left(
\sup_{\|z\|\leq L}|H_n(z)-G_n(z)|>\varepsilon \right) \leq
\exp\left\{ -C\left(\frac M2\right)^{\beta+1} \right\}.
\]
Letting $M\to\infty$ proves that $\sup_{\|z\|\leq L}|H_n(z)-G_n(z)|
\stk 0$.
Finally,
\[
\left|
\sup_{\|z\|\leq L}H_n(z)
-
\sup_{\|z\|\leq L}G_n(z)
\right|
\leq
\sup_{\|z\|\leq L}|H_n(z)-G_n(z)|.
\]
\eqref{eq:firstass} and  Slutsky's lemma therefore yield
\[
\sup_{\|z\|\leq L}H_n(z)
\xrightarrow{\mathcal D}
\Phi_L(N).
\]
This completes the proof.
\end{proof}


\subsection{The Poisson functional limit theorem}

We are now in a position to remove the restriction on the centre
variable and derive the full limit theorem for the smallest enclosing
radius. Lemma~\ref{lemma-7-8} gives convergence of the exact variational problem
on every fixed compact set, while Lemma~\ref{lem:weibull-centre-tightness} shows that the
optimization may, with arbitrarily large probability, be restricted
to a sufficiently large compact set. On the limiting side,
Lemma~\ref{lemma-7-7} allows us to remove the corresponding restriction.

\medskip

\begin{theorem}\label{hauptsatz-weibull}
Suppose that \eqref{eq:boundary-tail} holds. Then
\[
n^{1/(\beta + 1)}(1-R_n) \vertk \Phi(N),
\]
where
\[
\Phi(N) = \sup_{z\in\mathbb R^d} \inf_{(u,t)\in N}(t+z\cdot u),
\]
and $N$ is the Poisson point process from Theorem~\ref{thm:boundary-pp}.
\end{theorem}

\smallskip

\begin{proof} By Lemma~\ref{lemma-7-8}, $\sup_{\|z\|\leq L}H_n(z) \vertk 
\Phi_L(N)$ for every fixed $L>0$. Moreover, since the balls
$\{z:\|z\|\leq L\}$ increase to $\mathbb R^d$,
\[
\Phi_L(N) = \sup_{\|z\|\leq L}G(z) \uparrow \sup_{z\in\mathbb R^d}G(z) = \Phi(N) 
\]
almost surely as $L\to\infty$. By Lemma~\ref{lemma-7-7}, $\Phi(N)$ is finite
almost surely. On the other hand, Lemma~\ref{lem:weibull-centre-tightness} implies that, for every $\varepsilon>0$, there exists $L>0$ such that, for all sufficiently
large $n$,
\[
\mathbb P\left( \sup_{z\in\mathbb R^d}H_n(z) = \sup_{\|z\|\leq L}H_n(z) \right)
\geq1-\varepsilon.
\]
Hence the compactly restricted variables $\sup_{\|z\|\leq L}H_n(z)$ approximate
$\sup_{z\in\mathbb R^d}H_n(z)$ uniformly in probability for large
$n$, while $\Phi_L(N)$ converges almost surely to $\Phi(N)$ as
$L\to\infty$. The converging-together theorem therefore yields
$\sup_{z\in\mathbb R^d}H_n(z) \vertk \Phi(N)$.
The assertion now follows from the exact 
variational identity \eqref{eq:exact-weibull-variational}.
\end{proof}

\medskip

\begin{corollary}
Let $X_1,X_2,\ldots$ be i.i.d. uniformly distributed on $\mathbb B^d$. Then
\[
n(1-R_n)\vertk \Phi(N),
\]
where $N$ is the Poisson point process from Theorem~\ref{thm:boundary-pp}.
\end{corollary}
\smallskip

\begin{proof}
The assertion follows from Theorem~\ref{hauptsatz-weibull} and Remark~\ref{rem:uniform-ball}, since
\eqref{eq:boundary-tail} holds with $\beta=0$ and $C=d$.
\end{proof}

\medskip

\begin{remark}\label{rem:nondegenerate}
The distribution of $\Phi(N)$ has unbounded support and is therefore
nondegenerate. Indeed, putting
$T_0:=\inf\{t:(u,t)\in N\}$,
we have
\[
\mathbb{P}(T_0>s)=\exp\{-Cs^{\beta+1}\}, \qquad s>0.
\]
Thus, $T_0$ has unbounded support. Moreover, $\Phi(N)\geq G(0)=T_0$ almost surely, and hence
$\mathbb P(\Phi(N)>s)\geq\mathbb P(T_0>s) =\exp\{-Cs^{\beta+1}\}>0$ for every $s>0$.
Consequently, $\Phi(N)$ has unbounded support and is therefore
nondegenerate.
\end{remark}

\medskip

\begin{remark}\label{rem:scaling}
The constant $C$ in the boundary condition \eqref{eq:boundary-tail} can be
absorbed into the normalization. To see this, let $N_C$ denote a
Poisson point process on $\mathbb S^{d-1}\times[0,\infty)$ with intensity measure $C(\beta+1)t^\beta\,\D t\,\sigma_{d-1}(\D u)$,
and let $N_1$ denote the corresponding process with $C=1$.
By the scaling property of Poisson point processes, 
\[ 
N_C \stackrel{\mathcal D}{=} \left\{ \left(u,C^{-1/(\beta+1)}t\right):(u,t)\in N_1
\right\}.
\]
Consequently,
\[
\begin{split}
\Phi(N_C)
&\stackrel{\mathcal D}{=}
\sup_{z\in\mathbb R^d}
\inf_{(u,t)\in N_1}
\left(C^{-1/(\beta+1)}t+z\cdot u\right)\\
&=
C^{-1/(\beta+1)}
\sup_{y\in\mathbb R^d}
\inf_{(u,t)\in N_1}
(t+y\cdot u)\\
&=
C^{-1/(\beta+1)}\Phi(N_1),
\end{split}
\]
where in the second equality we put
$z=C^{-1/(\beta+1)}y$.
Thus Theorem~\ref{hauptsatz-weibull} may equivalently be written as
\[
(Cn)^{1/(\beta+1)}(1-R_n) \vertk \Phi(N_1).
\]
Hence the constant $C$ affects only the scale of the limiting distribution.
\end{remark}

\medskip

\begin{remark}\label{rem:compareweib}
It is instructive to compare the normalization in Theorem~\ref{hauptsatz-weibull}
with the corresponding normalization for the diameter. In the same
bounded-support setting, Lao (see Corollary~3.7 of \cite{lao10})
showed that
\[
n^{4/(4\beta+d+3)}(2-D_n) \vertk W,
\]
where $W$ has a nondegenerate Weibull distribution. Since
\[
\frac{1}{\beta+1}-\frac{4}{4\beta+d+3} = \frac{d-1}{(\beta+1)(4\beta+d+3)} >0,
\]
the deficit $1-R_n$ converges to zero on a strictly faster scale than $2-D_n$.
More precisely, Theorem~\ref{hauptsatz-weibull} implies
$n^{1/(\beta+1)}(1-R_n)=O_{\mathbb P}(1)$, and therefore
\[
n^{4/(4\beta+d+3)}(1-R_n) = n^{\,4/(4\beta+d+3)-1/(\beta+1)} \,n^{1/(\beta+1)}(1-R_n)
\stk 0.
\]
Since $2R_n-D_n = (2-D_n)-2(1-R_n)$, Slutsky's lemma yields
\begin{equation}\label{eq:finelimit}
n^{4/(4\beta+d+3)}(2R_n-D_n)
\vertk W.
\end{equation}
Thus, on the natural fluctuation scale of the diameter, the
difference between $2R_n$ and $D_n$ has the same limiting
distribution as the diameter deficit $2-D_n$.

This difference in scales reflects a geometric distinction between
the two extremal problems. The behaviour of $1-R_n$ is governed by
sample points lying close to the boundary whose directions form a
configuration capable of determining the enclosing ball. By
contrast, for $D_n$ to be close to $2$, two sample points must not
only lie close to the boundary but also have nearly antipodal
directions. This additional angular constraint leads to the slower
normalization of the diameter.
\end{remark}

\subsection{Numerical illustration}\label{subsecsimul}

Theorem~\ref{hauptsatz-weibull} gives the limiting distribution of the normalized
enclosing-radius deficit 
\[
a_n(1-R_n), \qquad a_n=n^{1/(\beta+1)},
\]
in terms of the Poisson functional $\Phi(N)$. Since the distribution
of this functional does not appear to admit a simple explicit
description, it is useful to examine the finite-sample behaviour of
the normalized statistic.

We therefore carried out a small simulation study in dimension
$d=2$ for several values of the boundary-decay parameter $\beta$.
For $\beta=0$, the observations were uniformly distributed on the
unit disk $\mathbb B^2$. Thus $C=2$ in 
\eqref{eq:boundary-tail}. For $\beta=1$ and $\beta=2$, the radial distribution was chosen such that
\[
\mathbb P(1-\rho\leq s)=s^{\beta+1},
\qquad 0\leq s\leq1,
\]
so that \eqref{eq:boundary-tail} holds with $C=1$. For each simulated sample we
computed the smallest enclosing radius $R_n$ and recorded the
normalized quantity
\[
a_n(1-R_n)=n^{1/(\beta+1)}(1-R_n).
\]

Table~\ref{tab:weibull-radius-sim} reports empirical means and selected quantiles for
$n=50,100,200$. For each fixed value of $\beta$, the entries are
already remarkably stable across the three sample sizes. In
particular, both the empirical means and the quantiles show only
minor variation with $n$. This provides a finite-sample illustration
of the convergence in distribution established in Theorem~\ref{hauptsatz-weibull}.

The numerical values corresponding to different values of $\beta$
should not, however, be compared directly. Apart from the change in
the boundary-decay exponent, the constant $C$ in \eqref{eq:boundary-tail} also
differs between the uniform model and the models used for
$\beta=1,2$. As shown in Remark~\ref{rem:scaling}, this constant changes the scale
of the limiting distribution. The relevant feature of Table~\ref{tab:weibull-radius-sim} is
therefore the stability with increasing $n$ for each fixed model,
rather than the relative size of the entries across different values
of $\beta$.

\begin{table}[ht]
\centering
\begin{tabular}{c|c|c|c|c|c|c}
$\beta$ & $n$ & mean & $q_{0.10}$ & $q_{0.25}$ & median & $q_{0.90}$ \\
\hline
0 & 50   & 1.372 & 0.480 & 0.771 & 1.218 & 2.467 \\
0 & 100  & 1.383 & 0.482 & 0.774 & 1.223 & 2.494 \\
0 & 200  & 1.388 & 0.483 & 0.776 & 1.227 & 2.505 \\
\hline
1 & 50   & 1.420 & 0.855 & 1.100 & 1.397 & 2.014 \\
1 & 100  & 1.440 & 0.865 & 1.112 & 1.414 & 2.047 \\
1 & 200  & 1.453 & 0.873 & 1.121 & 1.426 & 2.069 \\
\hline
2 & 50   & 1.215 & 0.863 & 1.028 & 1.214 & 1.568 \\
2 & 100  & 1.226 & 0.870 & 1.036 & 1.224 & 1.584 \\
2 & 200  & 1.234 & 0.875 & 1.042 & 1.232 & 1.594
\end{tabular}
\caption{Simulation results in dimension $d=2$ for the normalized quantity
$a_n(1-R_n)$, where $a_n=n^{1/(\beta+1)}$ ($10^6$ replications). For $\beta=0$ the sample is uniform
on $\mathbb B^2$. For $\beta=1,2$, the boundary depth $1-\rho$ satisfies
$\mathbb{P}(1-\rho\le s)=s^{\beta+1}$, $0\le s\le 1$.}
\label{tab:weibull-radius-sim}
\end{table}

Taken together, the simulation results indicate that the asymptotic
regime described by Theorem~~\ref{hauptsatz-weibull} is already visible for rather
moderate sample sizes in the planar examples considered here.



\section{The ratio $R_n/(D_n/2)$}\label{secratio}


The enclosing radius $R_n$ and the diameter $D_n$ satisfy the deterministic inequality
$D_n/2 \le R_n$.
It is therefore natural to ask under which conditions the two quantities are asymptotically equivalent, that is, whether $R_n/(D_n/2) \to 1$.
As shown in Remark~\ref{remcombine}, this convergence holds in probability in the planar
 Gaussian case. 
We now investigate the ratio $R_n/(D_n/2)$ under more general assumptions and first
consider the case of bounded support.

\begin{proposition}
Suppose that the distribution of $X$ has support equal to the closed unit ball
$\mathbb{B}^d = \{x \in \mathbb{R}^d : \|x\| \le 1\}$.
Then $R_n \to 1$ a.s., $D_n \to 2$ a.s., and therefore
\[
\frac{R_n}{D_n/2} \to 1 \quad \text{a.s.}
\]
\end{proposition}

\begin{proof}
Since all observations lie in $\mathbb{B}^d$, we have $R_n \le 1$ and $D_n \le 2$. 
Because the distribution has full support on $\mathbb{B}^d$, each neighbourhood of a boundary point of the
sphere $\mathbb{S}^{d-1}$ has positive probability. By independence, with probability one the sample
eventually contains points arbitrarily close to any prescribed boundary points.

Fix $\varepsilon > 0$. Choosing two antipodal points on the sphere shows that $D_n \ge 2 - \varepsilon$ eventually,
hence $D_n \to 2$ almost surely. 
Similarly, choosing the vertices of a regular simplex on
$\mathbb S^{d-1}$ and using continuity of the circumradius under
perturbations yields $R_n\ge 1-\varepsilon$ eventually,
so $R_n \to 1$ almost surely. The ratio statement follows.
\end{proof}

\medskip

The behaviour in the Fréchet regime is qualitatively different
from that in the bounded-support case.
As discussed in Section~\ref{secfrechet}, after
normalization the extreme observations converge to a Poisson point process, and both $R_n$ and
$D_n$ admit limit descriptions in terms of this limiting cloud. The following proposition shows
that the ratio $R_n/(D_n/2)$ converges in distribution to a functional of the limiting Poisson
configuration.

\medskip

\begin{proposition} Under the assumptions and with the notations of 
Theorem~\ref{thm:Frechet-radius},
\[
\frac{R_n}{D_n/2}
\vertk
\frac{\mathcal{C}(\{P_k: k \ge 1 \})}{\frac{1}{2}\sup_{1 \le i < j < \infty}\|P_i-P_j\|}.
\]
\end{proposition}

\smallskip

\noindent The limit is nondegenerate and, in particular, not concentrated at $1$.
To see this, fix $a>0$ and choose
\[
0<r_0<a\left(\sqrt{\frac{2(d+1)}{d}}-1\right).
\]
Since the intensity measure of the limiting Poisson process is finite on
$\{x:\|x\|>r_0\}$, events prescribing finitely many points in pairwise
disjoint Borel sets outside $B(0,r_0)$ and no further points outside
$B(0,r_0)$ have strictly positive probability.

First choose two sufficiently small disjoint neighbourhoods of two
antipodal points of norm $a$. On the event that each neighbourhood
contains exactly one Poisson point and that there are no further points
outside $B(0,r_0)$, all remaining Poisson points lie in $B(0,r_0)$.
For the unperturbed configuration consisting of the two antipodal points
together with $B(0,r_0)$, the circumradius equals $a$ and the diameter
equals $2a$. By continuity of the circumradius and diameter under small
perturbations of the two outer points, the ratio
\[
\frac{C(\{P_k:k\geq 1\})}
{\frac12\sup_{i<j}\|P_i-P_j\|}
\]
can therefore be made arbitrarily close to $1$ by choosing the two
neighbourhoods sufficiently small.

Next choose $d+1$ sufficiently small disjoint neighbourhoods of the
vertices of a regular simplex centred at the origin and inscribed in the
sphere of radius $a$. On the event that each neighbourhood contains
exactly one Poisson point and that there are no further points outside
$B(0,r_0)$, all remaining Poisson points again lie in $B(0,r_0)$.
For the unperturbed configuration consisting of the vertices of the
regular simplex together with $B(0,r_0)$, the circumradius equals $a$,
while the diameter equals $a\sqrt{2(d+1)/d}$.
Hence, by continuity, the above ratio can be made arbitrarily close to
$\sqrt{2d/(d+1)}>1$.
Choosing the neighbourhoods sufficiently small so that the two resulting
intervals for the ratio are disjoint, we obtain two disjoint events of
strictly positive probability on which the limiting ratio takes values
in these respective intervals. Hence the limiting random variable is
nondegenerate.


\smallskip

\begin{proof}
For $\varepsilon>0$, put
\[
Y_n^{(\varepsilon)} := C\bigl(\operatorname{conv}(\{0\}\cup N_n^{(\varepsilon)})\bigr),
\qquad
Z_n^{(\varepsilon)} := \operatorname{diam}\bigl(\{0\}\cup N_n^{(\varepsilon)}\bigr),
\]
and define analogously
\[
Y^{(\varepsilon)} := C\bigl(\operatorname{conv}(\{0\}\cup N^{(\varepsilon)})\bigr),
\qquad Z^{(\varepsilon)} := \operatorname{diam}\bigl(\{0\}\cup N^{(\varepsilon)}\bigr).
\]
For every fixed $\varepsilon>0$, Lemma~\ref{lem:Frechet-truncation} yields
$N_n^{(\varepsilon)} \vertk N^{(\varepsilon)}$.
Since the limiting truncated point process contains only finitely many
points almost surely, convergence of the truncated point processes
implies convergence of the corresponding finite point configurations
in the Hausdorff metric. Both the circumradius and the diameter are
continuous with respect to Hausdorff convergence of compact sets.
Consequently, by the continuous mapping theorem,
\[
\bigl(Y_n^{(\varepsilon)},Z_n^{(\varepsilon)}\bigr)
\xrightarrow{\mathcal D}
\bigl(Y^{(\varepsilon)},Z^{(\varepsilon)}\bigr)
\qquad\text{as }n\to\infty .
\]
It remains to remove the truncation. For the circumradius this has
already been done in Lemma~\ref{lem:Frechet-truncation}. Thus, for every $\delta>0$,
\[
\lim_{\varepsilon\downarrow0}
\limsup_{n\to\infty}
\mathbb P\left(
\left|c_n^{-1}R_n-Y_n^{(\varepsilon)}\right|>\delta
\right)
=0.
\]
We now establish the analogous statement for the diameter. With the
notation used in the proof of Lemma~\ref{lem:Frechet-truncation}, let
$S_n^{(\varepsilon)} = \{c_n^{-1}X_1,\ldots,c_n^{-1}X_n\}$
and $T_n^{(\varepsilon)} = \{0\}\cup \{c_n^{-1}X_k:\|X_k\|>\varepsilon c_n,\ 1\leq k\leq n\}$.
On the event $B_n^{(\varepsilon)}\neq\varnothing$, where
$B_n^{(\varepsilon)} = \{c_n^{-1}X_k:\|X_k\|\leq\varepsilon c_n,\ 1\leq k\leq n\}$,
we have
$d_H\bigl(S_n^{(\varepsilon)},T_n^{(\varepsilon)}\bigr)
\leq\varepsilon$.
Since the diameter functional satisfies
$\bigl|\operatorname{diam}(K)-\operatorname{diam}(L)\bigr|
\leq 2d_H(K,L)$ 
for non-empty compact sets $K,L$, it follows that, on
$\{B_n^{(\varepsilon)}\neq\varnothing\}$,
$\big|c_n^{-1}D_n-Z_n^{(\varepsilon)}\big| \leq 2\varepsilon$.
Moreover, as shown in the proof of Lemma~\ref{lem:Frechet-truncation},
$\mathbb P\bigl(B_n^{(\varepsilon)}=\varnothing\bigr)\to 0$
as $n\to\infty$
for every fixed $\varepsilon>0$. Hence, for every $\delta>0$,
\[
\lim_{\varepsilon\downarrow0} \limsup_{n\to\infty} 
\mathbb P\left(\left|c_n^{-1}D_n-Z_n^{(\varepsilon)}\right|>\delta \right) =0.
\]

For the limiting Poisson cloud, put
$S:=\{0\}\cup\{P_k:k\geq1\}$ and $T^{(\varepsilon)} := \{0\}\cup\{P_k:\|P_k\|>\varepsilon\}$.
Then
$d_H\bigl(S,T^{(\varepsilon)}\bigr)\leq\varepsilon$ a.s.
and therefore
$\left|
\operatorname{diam}(S)
- Z^{(\varepsilon)}  \right| \leq2\varepsilon$ a.s.

Together with the corresponding assertion for the circumradius from
Lemma~\ref{lem:Frechet-truncation}, this gives
\[
\bigl(Y^{(\varepsilon)},Z^{(\varepsilon)}\bigr)
\xrightarrow{\mathbb P}
\left(
C(\{P_k:k\geq1\}),
\sup_{1\leq i<j<\infty}\|P_i-P_j\|
\right)
\]
as $\varepsilon\downarrow0$. Here we have used that adjoining the
origin does not change either the circumradius or the diameter almost
surely, because the origin belongs almost surely to the closure of
$\{P_k:k\geq1\}$.

Hence, by the converging-together theorem applied in $\mathbb R^2$,
\[
\left(c_n^{-1}R_n,c_n^{-1}D_n\right)
\xrightarrow{\mathcal D}
\left(
C(\{P_k:k\geq1\}),
\sup_{1\leq i<j<\infty}\|P_i-P_j\|
\right).
\]
Since $\sup_{1\leq i<j<\infty}\|P_i-P_j\|>0$ almost surely, 
the mapping
$(x,y)\longmapsto x/(y/2)$ is continuous at the limiting random vector almost surely. A final application of the continuous mapping theorem therefore yields
\[
\frac{R_n}{D_n/2}
\vertk 
\frac{C(\{P_k:k\geq1\})}{\frac12\sup_{1\leq i<j<\infty}\|P_i-P_j\|},
\]
which proves the assertion.
\end{proof}

\medskip

\noindent We now turn to the Gumbel regime studied by Jammalamadaka and Janson \cite{jamsva15}.
In that setting, there exist sequences $a_n \to \infty$ and $b_n = o(a_n)$ such that
\[
\mathbb{P}(\|X\| > a_n + t b_n) = \frac{1+o(1)}{n} \E^{-t}
\]
uniformly for $t$ in a suitable range. In particular,
$(M_n - a_n)/b_n \vertk V$,
where $M_n = \max_{1 \le j \le n} \|X_j\|$ and $V$ has the Gumbel distribution.
Moreover, by Remark~1.3 in \cite{jamsva15},
\begin{equation}\label{jjconv}
\frac{2M_n - D_n}{b_n \log(a_n/b_n)} \xrightarrow{\mathbb{P}} \frac{d-1}{2}.
\end{equation}
Since $D_n \le 2R_n \le 2M_n$,
it follows that $0 \le 2R_n - D_n \le 2M_n - D_n$,
and hence the difference $2R_n - D_n$ is at most of order $b_n \log(a_n/b_n)$ in probability.

These observations immediately imply asymptotic equivalence of the
enclosing radius and half the diameter in the Gumbel regime.

\begin{proposition}
Let $X,X_1,X_2,\ldots$ be i.i.d.\ spherically symmetric random vectors
in $\mathbb{R}^d$ whose radial distribution belongs to the Gumbel
max-domain of attraction and satisfies the regularity assumptions of
Jammalamadaka and Janson \cite{jamsva15}. Then
\[
\frac{R_n}{D_n/2}\xrightarrow{\mathbb P}1.
\]
\end{proposition}

\begin{proof}
We  start with  \eqref{jjconv} and recall  $a_n\to\infty$ as well as $b_n=o(a_n)$. Since $(b_n/a_n)\log (a_n/b_n)\rightarrow 0$,
it follows that $2M_n-D_n=o_{\mathbb P}(a_n)$.
Moreover, $M_n/a_n\stk 1$, and hence
$D_n/a_n \stk 2$.
Since $0\leq 2R_n-D_n\leq 2M_n-D_n$,
we obtain
\[
0\leq
\frac{R_n}{D_n/2}-1
=
\frac{2R_n-D_n}{D_n}
\leq
\frac{2M_n-D_n}{D_n}
\xrightarrow{\mathbb P}0.
\]
This proves the assertion.
\end{proof}

The preceding result settles the first-order asymptotic relation
between $R_n$ and $D_n$ in the Gumbel regime. The finer behaviour of
their difference remains open. In the Weibull regime, the corresponding finer comparison is resolved by Remark~\ref{rem:compareweib}, which gives the limiting distribution of
$2R_n-D_n$ on the natural diameter scale.

\begin{problem}\label{prob:morma}
Under the assumptions of the preceding proposition, does there exist a
normalization under which $2R_n-D_n$
converges in distribution to a nondegenerate limit? In particular,
does
\[
\frac{2R_n-D_n}{b_n\log(a_n/b_n)}
\]
converge in distribution, and if so, what is the limiting law?
\end{problem}

The preceding results show that the asymptotic relationship between
$R_n$ and $D_n$ depends qualitatively on the extreme-value domain of
attraction of the radial distribution. In the bounded-support case,
$R_n/(D_n/2)\to1$ almost surely, and in the Weibull regime the finer
difference $2R_n-D_n$ admits the limit law \eqref{eq:finelimit}. In the
Fr\'echet case, the ratio converges to a nondegenerate random limit,
whereas in the Gumbel case it converges to $1$ in probability, with
the finer behaviour of $2R_n-D_n$ remaining open.

\section{Supporting configurations and asymptotic dominance}\label{secsuffcond}

A fundamental geometric feature of the smallest enclosing ball of a finite point set in
$\mathbb{R}^d$ is that it is supported by at most $d+1$ points. Recall that $K_n$ denotes
the number of sample  points supporting the smallest enclosing
ball of $X_1,\dots,X_n$.

The purpose of this section is to investigate the asymptotic behaviour of $K_n$.
In the planar case $d=2$, we establish rigorous results showing that the
three-point configuration becomes dominant as the sample size increases, in the sense that
$\mathbb{P}(K_n=3)\to 1$ as $n\to\infty$.
We then complement these results by a simulation study in dimensions $d=3$ and $d=4$,
which provides evidence that, more generally, the $(d+1)$-point configuration is
asymptotically dominant.

These findings motivate the conjecture that, under suitable regularity assumptions,
\begin{equation}\label{eq-d+1-dom}
\mathbb{P}(K_n=d+1)\to 1 \qquad \text{as } n\to\infty,
\end{equation}
for arbitrary dimension $d\ge 2$.

\subsection{The bounded support case}


In this subsection we assume that the support of the distribution of $X$ is the full unit disk $\mathbb{B}^2$. The following result shows that, under suitable assumptions on the radial distribution near the boundary of $\mathbb{B}^2$, the three-point configuration prevails with probability tending to one.

\medskip
\begin{proposition}\label{propkn2}
Let $X, X_1,X_2,\dots$ be i.i.d.\ random vectors in $\mathbb{R}^2$ with an absolutely continuous, spherically symmetric distribution supported on the closed unit disk $\mathbb{B}^2$.

Assume that there exist constants $\alpha>0$, $t_0\in(0,1)$ and $c_-,c_+>0$ such that
\[
c_- t^\alpha \le \mathbb{P}(\|X\|\ge 1-t)\le c_+ t^\alpha, \qquad 0<t\le t_0 .
\]
Then
\[
\lim_{n\to\infty}\mathbb{P}(K_n=3)=1.
\]
\end{proposition}

\smallskip
\begin{proof}
For $x,y\in\mathbb{B}^2$, let $D(x,y)$ denote the closed disk having the line segment $[x,y]$ as diameter, and define $\Delta(x,y):=1-\|x-y\|/2$.
Thus $D(x,y)$ has radius $1-\Delta(x,y)$. If $K_n=2$, the minimum enclosing circle of
$\{X_1,\dots,X_n\}$ equals $D(X_i,X_j)$ for some $1\le i<j\le n$, and therefore
$\{X_1,\dots,X_n\}\subset D(X_i,X_j)$.

Define $E_{ij}:=\{X_k\in D(X_i,X_j)\text{ for all }k\neq i,j\}$.
Then $\{K_n=2\}\subset\bigcup_{1\le i<j\le n}E_{ij}$, and
\begin{equation}\label{pkn2ob}
\mathbb{P}(K_n=2)\le \binom{n}{2}\mathbb{P}(E_{12})
\end{equation}
by exchangeability. We show that the right hand side of \eqref{pkn2ob} tends to zero as $n \to \infty$. Since $K_n \in \{2,3\}$ almost surely, we are done. 
Conditioning on $(X_1,X_2)$, we have
\[
\mathbb{P}(E_{12})
= \int_{\mathbb{B}^2} \int_{\mathbb{B}^2} \mathbb{P}\bigl(X_3\in D(x,y)\bigr)^{n-2}\,
\mathbb{P}_{X_1}(\D x)\,\mathbb{P}_{X_2}(\D y),
\]
where $\mathbb{P}_{X_j}$ denotes the distribution of $X_j$. In what follows we derive a lower bound for
$\mathbb{P}(X_3 \notin D(x,y))$, which in turn yields an upper bound for $\mathbb{P}(X_3 \in D(x,y))$.
\end{proof}
\medskip
\noindent

\begin{lemma}\label{lemma1k}
There exists a constant $a_1>0$ such that for all $x,y\in \mathbb{B}^2$, writing
$\delta:=\Delta(x,y)=1-\|x-y\|/2 \in (0,1]$, we have
\[
\mathbb{P}(X_3 \notin D(x,y)) \ge a_1 \delta^\alpha.
\]
\end{lemma}

\smallskip
\noindent
\begin{proof}
Let $m=(x+y)/2$ be the centre of $D(x,y)$ and $\rho=1-\delta$ its radius. Write
$s=\|m\|$, let $v=m/s$ if $s>0$, and choose any unit vector $v$ otherwise. Fix $\eta\in(0,\pi/4)$,
set $h(\delta):=\min(t_0,\delta/2)$, and consider the set
$S(\delta,v):=\{r\E^{i\theta}: 1-h(\delta)\le r\le 1,\ |\theta-(\arg v+\pi)|\le \eta\}$, which is a boundary sector of the unit disk pointing in the direction opposite to $v$.

We claim that $S(\delta,v)\subset B(0,1)\setminus D(x,y)$ for all $\delta\in(0,1]$.

If $s=0$, then $m=0$ and $D(x,y)=B(0,\rho)$ with $\rho=1-\delta$. Hence for any
$z\in S(\delta,v)$ we have $\|z\|\ge 1-h(\delta)\ge 1-\delta/2>1-\delta=\rho$, so $z\notin D(x,y)$.

Now suppose that $s>0$, and let $z=r\E^{i\theta}\in S(\delta,v)$. Then the angle between $z$ and $m$
is at least $\pi-\eta$. By the cosine law,
\[
\|z-m\|^2 = r^2 + s^2 - 2rs\cos(\angle(z,m)) \ge r^2 + s^2 + 2rs\cos\eta.
\]
Since $r\ge 1-h(\delta)\ge 1-\delta/2$ and $\cos\eta>0$, it follows that
\[
\|z-m\|^2 \ge (1-\delta/2)^2 + s^2 + 2(1-\delta/2)s\cos\eta.
\]
Moreover,
\[
(1-\delta/2)^2 + s^2 - (1-\delta)^2 = \delta + s^2 - \tfrac{3}{4}\delta^2 > 0
\]
for all $\delta\in(0,1]$. Hence $\|z-m\|^2 > (1-\delta)^2=\rho^2$, so $z\notin D(x,y)$.
Thus $S(\delta,v)\subset B(0,1)\setminus D(x,y)$.

By spherical symmetry,
\[
\mathbb{P}(X_3 \in S(\delta,v))
= \frac{2\eta}{2\pi}\,\mathbb{P}(\|X_3\|\ge 1-h(\delta))
= \frac{\eta}{\pi}\,\mathbb{P}(\|X_3\|\ge 1-h(\delta)).
\]
Since $h(\delta)\le t_0$, the lower tail bound yields 
$\mathbb{P}(\|X_3\|\ge 1-h(\delta)) \ge c_- h(\delta)^\alpha$.
If $\delta\le 2t_0$, then $h(\delta)=\delta/2$, so $h(\delta)^\alpha = 2^{-\alpha}\delta^\alpha$.
If $\delta>2t_0$, then $h(\delta)=t_0$, and since $\delta\le 1$, we have
$h(\delta)^\alpha = t_0^\alpha \ge t_0^\alpha \delta^\alpha$.
Hence, for all $\delta\in(0,1]$,
\[
\mathbb{P}(X_3 \notin D(x,y))
\ge \frac{\eta}{\pi} c_- \min\{2^{-\alpha}, t_0^\alpha\}\,\delta^\alpha.
\]
This proves the assertion with $a_1 = (\eta/\pi)c_- \min\{2^{-\alpha}, t_0^\alpha\}$.
\end{proof}

\medskip
\noindent

\begin{lemma}\label{lemma2k}
There exists $a_2>0$ such that for all sufficiently small $\delta>0$,
\[
\mathbb{P}\bigl(\Delta(X_1,X_2)\le \delta\bigr) \le a_2 \delta^{2\alpha+1/2}.
\]
\end{lemma}

\smallskip
\noindent
\begin{proof}
Write $X_j=\rho_j \E^{i\Theta_j}$ ($j=1,2$), where $\Theta_1, \Theta_2, \rho_1, \rho_2$ are independent, and $\Theta_1, \Theta_2$ are i.i.d.\ uniform on $[0,2\pi)$.
Let $T_j=1-\rho_j$ ($j=1,2$). If $\Delta(X_1,X_2)\le \delta$, then $\|X_1-X_2\|\ge2(1-\delta)$,
while $\|X_1-X_2\|\le \rho_1+\rho_2$. Hence $T_1+T_2\le 2\delta$, so $\max(T_1,T_2) \le 2\delta$.
Let $\Psi=\pi-|\Theta_1-\Theta_2|_{2\pi}$, where $|\, \cdot \, |_{2\pi}$ denotes the circular distance modulo $2\pi$. Since $\min(\rho_1,\rho_2)\ge 1-2\delta$ and $\|X_1-X_2\|\ge 2(1-\delta)$, the cosine law gives
\[
2\rho_1\rho_2(1-\cos\Psi) = (\rho_1+\rho_2)^2-\|X_1-X_2\|^2
\le 4-4(1-\delta)^2 \le 8\delta.
\]
Hence $1-\cos\Psi \le 4\delta/(1-2\delta)^2$.
For sufficiently small $\delta$, this implies $\Psi\le 1$. Using $1-\cos\Psi\ge \Psi^2/4$, $0\le \Psi\le 1$, we obtain $ \Psi\le C\sqrt{\delta} $ for some constant $C>0$. Thus
$\{\Delta(X_1,X_2)\le \delta\} \subset \{T_1\le 2\delta,\ T_2\le 2\delta,\ \Psi\le C\sqrt{\delta}\}$.
By independence, $\mathbb{P}(\Delta(X_1,X_2)\le \delta)
\le \mathbb{P}(T_1\le 2\delta)\mathbb{P}(T_2\le 2\delta)
\mathbb{P}(\Psi\le C\sqrt{\delta})$, which, together with the upper tail bound and uniformity of the angles, yields
\[
\mathbb{P}(\Delta(X_1,X_2)\le \delta)\le a_2 \delta^{2\alpha+1/2}.
\]
\end{proof}

\medskip
\noindent

\begin{proof}[Proof of Proposition~\ref{propkn2}, continued]
Returning to the main argument, let
\[
\delta_n=\Bigl(\frac{A\log n}{n}\Bigr)^{1/\alpha},
\]
where $A>4/a_1$, and $a_1$ is the constant from Lemma~\ref{lemma1k}. Conditioning on
$(X_1,X_2)$, we have
\[
\mathbb{P}(E_{12}\mid X_1,X_2)=\mathbb{P}\bigl(X_3\in D(X_1,X_2)\bigr)^{\,n-2}.
\]
On the event $\{\Delta(X_1,X_2)\le \delta_n\}$, we use the trivial bound
$\mathbb{P}(E_{12}\mid X_1,X_2)\le 1$. On the event $\{\Delta(X_1,X_2)>\delta_n\}$, Lemma~\ref{lemma1k} yields
$\mathbb{P}(E_{12}\mid X_1,X_2)
\le \bigl(1-a_1\Delta(X_1,X_2)^\alpha\bigr)^{n-2}
\le \bigl(1-a_1\delta_n^\alpha\bigr)^{n-2}$.
Taking expectations gives
\[
\mathbb{P}(E_{12}) \le \mathbb{P}\bigl(\Delta(X_1,X_2)\le \delta_n\bigr)
+ \bigl(1-a_1\delta_n^\alpha\bigr)^{n-2}.
\]
Using $1-t\le \E^{-t}$ for $t\ge 0$, we obtain $\bigl(1-a_1\delta_n^\alpha\bigr)^{n-2}
\le \exp\{-a_1(n-2)\delta_n^\alpha\}$. By Lemma~\ref{lemma2k},
$\mathbb{P}\bigl(\Delta(X_1,X_2)\le \delta_n\bigr) \le
a_2\delta_n^{2\alpha+1/2}$ for all sufficiently large $n$. Since $\delta_n^\alpha=A\log n/n$, we have
$\exp\{-a_1(n-2)\delta_n^\alpha\}\le n^{-a_1A/2}$ for all sufficiently large $n$. Thus
$\mathbb{P}(E_{12})\le a_2\delta_n^{2\alpha+1/2}+n^{-a_1A/2}$.

Multiplying by $\binom{n}{2}$ yields
\begin{equation}\label{endproof51}
\mathbb{P}(K_n=2)\le \binom{n}{2}\bigl(a_2\delta_n^{2\alpha+1/2}+n^{-a_1A/2}\bigr).
\end{equation}
Since
\[
n^2\delta_n^{2\alpha+1/2}
= A^{2+1/(2\alpha)}(\log n)^{2+1/(2\alpha)}\,n^{-1/(2\alpha)}
\]
converges to $0$, and since $A>4/a_1$ ensures that $n^{2-a_1A/2}\to 0$, the right-hand side of \eqref{endproof51} tends to zero, which was to be shown.
\end{proof}

\subsection{The Gaussian case}


The bounded-support result of Proposition~\ref{propkn2} suggests that, in the planar case,
the three-point configuration may also dominate in other regimes. The next theorem confirms
this in the Gaussian case treated by Daniels~\cite{dan52}. In contrast to the geometric argument
used in the bounded-support case, the Gaussian case requires an analytic approach based on
Daniels' integral representation.

\begin{theorem}\label{thm:gaussian-Kn3}
Let $X,X_1,X_2,\dots$ be i.i.d.\ standard normal random vectors in $\mathbb{R}^2$. Then
\[
\lim_{n\to\infty}\mathbb{P}(K_n=3)=1.
\]
\end{theorem}

\smallskip
\noindent 

\begin{proof}
Write $R_n$ for the radius and $C_n$ for the centre of the minimum
enclosing circle. As in Daniels~\cite{dan52}, let
$ P(r,\rho):=\mathbb P\bigl(X_1\in B(c,r)\bigr)$, $\|c\|=\rho$,
denote the probability that a standard normal point falls into a disk
of radius $r$ whose centre has distance $\rho$ from the origin. By
spherical symmetry, this probability depends only on $r$ and $\rho$.

Daniels distinguishes two mutually exclusive cases: either the minimum
enclosing circle is determined by three sample points forming an acute
triangle, or it is determined by two diametrically opposite sample
points. We first derive the contribution of the latter case.

To this end, fix two labelled sample points. Write the centre in polar coordinates as
$c=(\rho\cos\phi,\rho\sin\phi)$, and parametrize two diametrically opposite points on the circle by
\[
 \begin{split}
 x_1&=c+r(\cos\theta,\sin\theta),\\
 x_2&=c-r(\cos\theta,\sin\theta),
 \end{split}
 \qquad
 0\leq\phi<2\pi,\quad 0\leq\theta<2\pi.
\]
Here the two supporting points are regarded as labelled. Hence
$\theta$ ranges over $[0,2\pi)$: replacing $\theta$ by
$\theta+\pi$ interchanges the two labelled points and therefore
corresponds to a different ordered configuration. The passage to
unordered supporting pairs is accounted for below by the factor
$\binom{n}{2}$.

The absolute value of the Jacobian of the transformation
$(\rho,\phi,r,\theta)\mapsto(x_1,x_2)$ is $4\rho r$, while
$\|x_1\|^2+\|x_2\|^2=2(\rho^2+r^2)$.
Since the joint density of $(X_1,X_2)$ is
\[
 (2\pi)^{-2}  \exp\left\{-\frac{\|x_1\|^2+\|x_2\|^2}{2}\right\},
\]
integration over $\phi,\theta\in[0,2\pi)$ gives, for a fixed labelled
pair, $4\rho r\,\exp(-(\rho^2+r^2))\,\D \rho\,\D r$.
The remaining $n-2$ sample points must all lie in $B(c,r)$, which
contributes the factor $P(r,\rho)^{n-2}$. Finally, there are
$\binom{n}{2}$ choices for the unordered pair of supporting sample
points. Consequently, the infinitesimal probability corresponding to
the two-point case is
\[
 \D F_{n,2}(r,\rho)
 =
 2n(n-1)\rho r\,{\rm e}^{-(\rho^2+r^2)}
 P(r,\rho)^{n-2}\, \D \rho\,\D r.                    
\]
Thus
\begin{equation}\label{eq:dan52kn2}
 \mathbb P(K_n=2)  =  2n(n-1)  \int_0^\infty\int_0^\infty
 \rho r\,{\rm e}^{-(\rho^2+r^2)}  P(r,\rho)^{n-2}\, \D \rho\,\D r.                    
\end{equation}
We split the integral according to whether the centre-distance $\rho$
is at most $1$ or larger than $1$. By Lemma~\ref{lem:tight-centres} and 
Remark~\ref{sec3:rem38}, we have $\|C_n\|=O_{\mathbb P}\bigl((\log n)^{-1/2}\bigr)$,
and hence $\mathbb P(\|C_n\|>1)\to0$. It therefore remains to show that
$ I_n:=\mathbb P(K_n=2,\|C_n\|\leq1)\rightarrow 0$.
By \eqref{eq:dan52kn2},
\[
 I_n  =  2n(n-1)  \int_0^1\int_0^\infty  \rho r\,{\rm e}^{-(\rho^2+r^2)}
 P(r,\rho)^{n-2}\, \D r\, \D \rho.
\]

We next derive a lower bound for $1-P(r,\rho)$. From the first equation
on page 138 of Daniels~\cite{dan52},
\[
 1-P(r,\rho)  =  {\rm e}^{-\rho^2/2}  \int_r^\infty  s\,{\rm e}^{-s^2/2}I_0(\rho s)\, \D s,
\]
where $I_0$ denotes the modified Bessel function of order zero. Since
$I_0$ is increasing on $[0,\infty)$, it follows that 
$1-P(r,\rho)  \ge \exp(-(\rho^2+r^2)/2)I_0(\rho r)$.
For $0\leq\rho\leq1$, this implies
$ 1-P(r,\rho)  \geq  {\rm e}^{-1/2}{\rm e}^{-r^2/2}I_0(\rho r)$. Using
\[
 I_0(u)
 =\sum_{k=0}^\infty\frac{(u^2/4)^k}{(k!)^2}  \geq1+\frac{u^2}{4}
 \geq\frac14(1+u^2),  \qquad u\geq0,
\]
and putting $c_0:={\rm e}^{-1/2}/4$, we conclude that
$ 1-P(r,\rho)  \geq  c_0 \exp(-r^2/2)(1+r^2\rho^2)$ for 
 $0\leq\rho\leq 1$ and $r\geq 0$.
Since $1-u\leq{\rm e}^{-u}$ for $u\geq0$, it follows that
\[
 P(r,\rho)^{n-2}  \leq  \exp\left\{ -c_0(n-2){\rm e}^{-r^2/2}(1+r^2\rho^2)  \right\}.
\]
Set $\lambda_n(r):=c_0(n-2)\exp(-r^2/2)$. Since ${\rm e}^{-\rho^2}\leq1$, we obtain
\[
 I_n
 \le  Cn^2\int_0^\infty  r{\rm e}^{-r^2}{\rm e}^{-\lambda_n(r)}  \left(
 \int_0^1  \rho{\rm e}^{-\lambda_n(r)r^2\rho^2}\, \D \rho  \right) \D r,
\]
where $C>0$ denotes an absolute constant that may vary from line to
line but does not depend on $n$.

We split the last integral at $r=1$ and write
$I_n\leq I_{n,1}+I_{n,2}$. For $0<r\leq1$, the inner integral is at
most $1/2$, while $\lambda_n(r)  \geq c_0(n-2){\rm e}^{-1/2}$.
Hence $ I_{n,1}  \leq  Cn^2 \exp(-c_1n)  \rightarrow 0$ for some constant $c_1>0$.
For $r\geq1$, the inner integral is bounded by
$1/(\lambda_n(r)r^2)$. Therefore
\[
 I_{n,2}  \leq   Cn\int_1^\infty  \frac{{\rm e}^{-r^2/2}}{r}
 {\rm e}^{-\lambda_n(r)}\, \D r.
\]
With the substitution $y=\lambda_n(r)=c_0(n-2)\exp(-r^2/2)$ we obtain
\[
 I_{n,2}  \leq  C\int_0^{c_0(n-2){\rm e}^{-1/2}}  \frac{{\rm e}^{-y}}
 {2\log\{c_0(n-2)/y\}}\, \D y.
\]
Split this integral at $y=\sqrt n$. For $0<y\leq\sqrt n$,
\[
 \log\frac{c_0(n-2)}{y}  \geq\frac12\log n+O(1),
\]
and hence
\[
 \int_0^{\sqrt n}  \frac{{\rm e}^{-y}}{2\log\{c_0(n-2)/y\}}\, \D y 
 \leq\frac{C}{\log n}.
\]
On the remaining interval,
$y\leq c_0(n-2){\rm e}^{-1/2}$ implies $c_0(n-2)/y\geq{\rm e}^{1/2}$, so that
\[
 \int_{\sqrt n}^{c_0(n-2){\rm e}^{-1/2}}
 \frac{{\rm e}^{-y}}  {2\log\{c_0(n-2)/y\}}\, \D y
 \leq  \int_{\sqrt n}^\infty{\rm e}^{-y}\,dy  =  {\rm e}^{-\sqrt n}.
\]
Consequently,
\[
 I_{n,2}  \leq  \frac{C}{\log n}+C{\rm e}^{-\sqrt n}  \rightarrow 0.
\]
Together with $I_{n,1}\to 0$, this gives $I_n\to0$. Since
$\mathbb P(\|C_n\|>1)\to0$, we conclude that
$\mathbb P(K_n=2)\to0$. As $K_n\in\{2,3\}$ almost surely,
we have $\mathbb P(K_n=3)\rightarrow 1$, which proves the assertion.
\end{proof}

\medskip

\begin{remark}
The representation \eqref{eq:dan52kn2} agrees with the infinitesimal probability
formula (3.5) of Daniels~\cite{dan52}. In Daniels' notation,
$\D Q_2$ is the probability element corresponding to two labelled
points in diametrically opposite positions on the circle. For unit
variance, his formula (3.6) gives
\[
\D Q_2=4{\rm e}^{-(r^2+\rho^2)}r\rho\, \D r\,\D \rho,
\]
and multiplication by $\binom{n}{2}P(r,\rho)^{n-2}$ yields \eqref{eq:dan52kn2}.
\end{remark}


\subsection{Simulation evidence in higher dimensions and a conjecture}\label{subsechigh}


The preceding results establish asymptotic dominance of the three-point configuration in the planar case under assumptions that ensure a strong concentration of extreme points near the boundary of the support. It is natural to ask whether an analogous phenomenon holds in higher dimensions, both for bounded-support distributions and for distributions in the Gumbel max-domain of attraction.

In contrast, for distributions in the Fréchet max-domain of attraction, where extreme points occur as isolated outliers far from the bulk of the sample, a different geometric behaviour is to be expected. This is consistent with the simulation results reported in Tables~\ref{tab:frechet-support} and~\ref{tabsec3-2}, which do not provide evidence for dominance of the $(d+1)$-point configuration in this case.

To explore this question, we carried out a simulation study in dimensions $d=3$ and $d=4$
for both bounded-support and Gaussian models. For the bounded-support case, we generated
i.i.d.\ random vectors uniformly distributed on the unit ball $\mathbb{B}^d$, while for the
Gaussian case we generated i.i.d.\ standard normal random vectors in $\mathbb{R}^d$.
In each case, we computed the smallest enclosing ball and recorded the number $K_n$
of supporting points. Each entry in the tables is based on $10^6$ replications.

The resulting estimates are shown in Tables~\ref{tab:K3-unif-gauss} and~\ref{tab:K4-unif-gauss}. Although the sample
sizes are moderate, the results reveal a clear pattern in the
bounded-support case: in dimension $d=3$, the probability of the four-point
configuration increases markedly with $n$, and in dimension $d=4$ the
five-point configuration shows the same pronounced tendency. For the
Gaussian model, the probability of the $(d+1)$-point configuration also
increases with $n$, but much more slowly. This effect is particularly
pronounced in dimension $d=4$, where the four-point configuration is still
more frequent than the five-point configuration at $n=200$.

\begin{table}[!ht]
\centering
\begin{tabular}{llcccc}
Model & $n$ & $\widehat{\mathbb{P}}(K_n=2)$ & $\widehat{\mathbb{P}}(K_n=3)$ &
$\widehat{\mathbb{P}}(K_n=4)$ \\
\hline
Uniform on $\mathbb B^3$ & $20$  & $0.027$ & $0.318$ & $0.656$ \\
Uniform on $\mathbb B^3$ & $50$  & $0.005$ & $0.155$ & $0.840$ \\
Uniform on $\mathbb B^3$ & $100$ & $0.001$ & $0.083$ & $0.915$ \\
Uniform on $\mathbb B^3$ & $200$ & $0.000$ & $0.044$ & $0.956$ \\
\hline
Gaussian in $\mathbb{R}^3$ & $20$  & $0.141$ & $0.515$ & $0.344$ \\
Gaussian in $\mathbb{R}^3$ & $50$  & $0.094$ & $0.459$ & $0.447$ \\
Gaussian in $\mathbb{R}^3$ & $100$ & $0.072$ & $0.422$ & $0.505$ \\
Gaussian in $\mathbb{R}^3$ & $200$ & $0.058$ & $0.390$ & $0.552$ \\
\hline
\end{tabular}
\smallskip
\caption{Estimated probabilities $\widehat{\mathbb{P}}(K_n=j)$ for the support size of the smallest
enclosing ball in dimension $d=3$ ($10^6$ replications).}
\label{tab:K3-unif-gauss}
\end{table}

\begin{table}[!ht]
\centering
\begin{tabular}{llccccc}
Model & $n$ & $\widehat{\mathbb{P}}(K_n=2)$ & $\widehat{\mathbb{P}}(K_n=3)$ &
$\widehat{\mathbb{P}}(K_n=4)$ & $\widehat{\mathbb{P}}(K_n=5)$ \\
\hline
Uniform on $\mathbb B^4$ & $20$  & $0.005$ & $0.099$ & $0.432$ & $0.465$ \\
Uniform on $\mathbb B^4$ & $50$  & $0.000$ & $0.022$ & $0.249$ & $0.728$ \\
Uniform on $\mathbb B^4$ & $100$ & $0.000$ & $0.006$ & $0.144$ & $0.850$ \\
Uniform on $\mathbb B^4$ & $200$ & $0.000$ & $0.002$ & $0.077$ & $0.921$ \\
\hline
Gaussian in $\mathbb{R}^4$ & $20$  & $0.065$ & $0.349$ & $0.440$ & $0.146$ \\
Gaussian in $\mathbb{R}^4$ & $50$  & $0.037$ & $0.260$ & $0.468$ & $0.235$ \\
Gaussian in $\mathbb{R}^4$ & $100$ & $0.026$ & $0.214$ & $0.466$ & $0.293$ \\
Gaussian in $\mathbb{R}^4$ & $200$ & $0.019$ & $0.180$ & $0.458$ & $0.342$ \\
\hline
\end{tabular}
\smallskip
\caption{Estimated probabilities $\widehat{\mathbb{P}}(K_n=j)$ for the support size of the smallest
enclosing ball in dimension $d=4$ ($10^6$ replications)}.
\label{tab:K4-unif-gauss}
\end{table}

Taken together, the simulations reveal a clear difference between the
bounded-support and Gaussian settings. In the bounded-support case, the
$(d+1)$-point configuration becomes dominant already for moderate sample
sizes in dimensions $d=3$ and $d=4$, and the numerical results suggest a
comparatively rapid approach to $(d+1)$-point dominance. By contrast, for
the Gaussian model, although the estimated probability of the
$(d+1)$-point configuration increases consistently with $n$, this tendency
is much slower, and competing configurations retain substantial
probabilities even for $n=200$. This is especially apparent in dimension
$d=4$, where the four-point configuration is still more frequent than the
five-point configuration at the largest sample size considered.

These findings support the view that, in the bounded-support case, the
asymptotic geometry of the smallest enclosing ball is governed by points
accumulating near the boundary of the support, leading to dominance of the
$(d+1)$-point configuration.

This motivates the following conjecture.

\begin{conjecture}\label{conj:dplus1-bounded}
Let $X,X_1,X_2,\dots$ be i.i.d.\ absolutely continuous spherically symmetric random vectors in
$\mathbb{R}^d$, $d\ge 2$.
Suppose that the distribution of $X$ is supported on the closed unit ball $\mathbb B^d$ and that
there exist constants $\alpha>0$, $t_0\in(0,1)$ and $c_-,c_+>0$ such that
\[
c_- t^\alpha \le \mathbb{P}(\|X\|\ge 1-t)\le c_+ t^\alpha, \qquad 0<t\le t_0.
\]
Then \eqref{eq-d+1-dom} holds.
\end{conjecture}

The conjecture is also consistent with a simple geometric heuristic.
In the bounded-support setting, extreme observations concentrate near the boundary sphere.
A generic collection of $d+1$ such points determines a simplex-like configuration whose
circumscribed ball is close to the unit ball, whereas configurations involving fewer than
$d+1$ supporting points require additional geometric constraints and may therefore be
expected to occur with smaller probability as the sample size increases.

\noindent For distributions in the Gumbel max-domain of attraction, the situation is less clear.
Although Theorem~\ref{thm:gaussian-Kn3} establishes asymptotic dominance of the three-point configuration
in the planar Gaussian case, the higher-dimensional behaviour is not yet understood.

\begin{problem}
\label{prob:gumbel-support}
Let $X,X_1,X_2,\dots$ be i.i.d.\ spherically symmetric random vectors in $\mathbb{R}^d$
whose radial distribution belongs to the Gumbel max-domain of attraction and satisfies
the regularity assumptions of Jammalamadaka and Janson~\cite{jamsva15}.
Is it true that \eqref{eq-d+1-dom} holds?
\end{problem}


\section{Open problems and further directions}\label{secproblems}

The results of the preceding sections provide limit theorems for the
smallest enclosing radius in the three classical extreme-value
regimes, but several questions concerning the corresponding limit
laws and the underlying geometry remain open.

\medskip
\noindent
\textit{(1) Explicit characterization of the limiting Poisson
functionals.} In the Weibull case, Theorem~\ref{hauptsatz-weibull} identifies the limiting law of the normalized enclosing-radius deficit as the distribution of
\[
\Phi(N) = \sup_{z\in\mathbb R^d} \inf_{(u,t)\in N}(t+z\cdot u).
\]
Even in the planar uniform case, it remains an interesting problem to obtain a more explicit description of this distribution, for example in terms of its distribution function, density, moments, or tail behaviour. A related problem arises in the Gumbel regime: while the limiting law can be identified explicitly in dimension
$d=2$ by Daniels' formula, no corresponding explicit
characterization appears to be available in dimensions $d\geq3$.

\medskip
\noindent
\textit{(2) Limiting behaviour of the centre.}
The proofs in Sections~\ref{secgumbelkotz} and~\ref{secweib} establish tightness of the centre of the smallest enclosing ball on the natural fluctuation scale. In
particular, in the Weibull setting, $a_n C_n=O_{\mathbb P}(1)$.
It is natural to ask whether the rescaled centres themselves converge
in distribution. In the Weibull case, a plausible candidate would be
a maximizer of the random function
\[
z\mapsto \inf_{(u,t)\in N}(t+z\cdot u).
\]
This raises, in particular, the question whether this maximizer is
almost surely unique. Analogous questions may be asked in the
Gumbel regime, and joint limit theorems for the radius and the centre
would provide a more complete stochastic-geometric description of
the smallest enclosing ball.

\medskip
\noindent
\textit{(3) Supporting configurations in higher dimensions.}
Section~\ref{secsuffcond} shows that the three-point configuration is asymptotically
dominant in dimension $d=2$ for the bounded-support and Gaussian models considered there. The simulation results in dimensions $d=3$ and $d=4$ suggest that, in the bounded-support case, the $(d+1)$-point configuration may be asymptotically dominant more
generally. Proving Conjecture~\ref{conj:dplus1-bounded} remains open. It would also be of
interest to determine whether an analogous result holds in the
Gumbel regime, as asked in Problem~\ref{prob:gumbel-support}, and to relate the supporting
sample points to the active constraints of the limiting variational
problem.

\medskip
\noindent
\textit{(4) Finer comparison of the enclosing radius and the
diameter.}
The first-order relation between $R_n$ and $D_n$ is now well
understood in the three regimes considered here. In the Gumbel case,
however, the finer behaviour of
$2R_n-D_n$ remains open. In particular, Problem~\ref{prob:morma} asks whether there is a
natural normalization under which this difference converges in
distribution to a nondegenerate limit. Resolving this question would
clarify which geometric features distinguish the enclosing-radius
and diameter problems beyond first order.

\medskip
\noindent
\textit{(5) Beyond spherical symmetry.}
The point-process arguments developed in this paper rely strongly on
the spherical decomposition $X=\rho U$ with a uniform direction
$U$. It would be interesting to extend the theory to distributions
with non-uniform angular components and, more generally, to
distributions supported on smooth convex bodies other than the
Euclidean ball. In such settings both the local boundary geometry and
the angular intensity of extreme observations should enter the
limiting point process and the corresponding variational functional.

More broadly, it would be desirable to understand to what extent the
Poisson functionals arising in the Gumbel, Fr\'echet, and Weibull
regimes can be viewed within a unified stochastic-geometric
framework. Although the limiting point processes are quite different
in the three cases, the smallest enclosing ball is in each regime
determined asymptotically by a geometric optimization over an
extremal Poisson configuration.

\section*{Acknowledgements}

The authors made use of ChatGPT (OpenAI), accessed through the ChatGPT web interface on multiple occasions during the preparation of this manuscript in 2026. The tool was used interactively for language refinement, editorial assistance, and the discussion and formulation of mathematical arguments. Different model versions available through ChatGPT were used during this period. The code used for the simulations was generated using Claude (Anthropic). All mathematical statements, proofs, numerical results, references, and conclusions were reviewed by the authors, who take full responsibility for the content of the paper.

This work was supported by the Initiative and Networking Fund of the
Helmholtz Association under the call Helmholtz Young Investigator Groups
(VH-NG-19-34, DataMat). The authors gratefully acknowledge the
scientific support and HPC resources provided by the German Aerospace
Center (DLR). The HPC system CARO is partially funded by the Ministry of
Science and Culture of Lower Saxony and the Federal Ministry for Economic
Affairs and Climate Action.

\end{document}